\documentclass[11pt,a4paper]{amsart}
\usepackage[T1]{fontenc}
\usepackage{amsmath,amssymb,mathtools,mathrsfs}
\usepackage{array,booktabs,graphicx,adjustbox,xcolor}
\usepackage{enumerate}
\usepackage[hypcap=false,justification=raggedright,singlelinecheck=true]{caption}
\usepackage[numbers,square,sort&compress]{natbib}
\usepackage[a4paper,left=25mm,right=25mm,top=27mm,bottom=27mm,heightrounded]{geometry}
\usepackage{microtype}
\usepackage[colorlinks=true,citecolor=magenta!70!black,
  linkcolor=blue!65!black,urlcolor=violet!70!black]{hyperref}
\hypersetup{
  pdftitle={On minimal singular vectors of universal affine vertex algebras of types B and C},
  pdfauthor={Cuipo Jiang and Tong Wu},
  pdfsubject={Universal affine vertex algebras and minimal singular vectors},
  pdfkeywords={affine Lie algebra, singular vector, conformal weight, Jantzen filtration}
}
\newcommand{\CC}{\mathbb C}
\newcommand{\ZZ}{\mathbb Z}
\newcommand{\Stab}{\operatorname{Stab}}

\theoremstyle{plain}
\newtheorem{theorem}{Theorem}[section]
\newtheorem{lem}[theorem]{Lemma}
\newtheorem{prop}[theorem]{Proposition}
\newtheorem{cor}[theorem]{Corollary}
\makeatletter
\@addtoreset{theorem}{subsection}
\makeatother
\renewcommand{\thetheorem}{%
  \ifnum\value{subsection}=0
    \thesection.\arabic{theorem}%
  \else
    \thesubsection.\arabic{theorem}%
  \fi
}
\numberwithin{equation}{section}

\newcommand{\kg}{\mathfrak{g}}

\def\leq{\leqslant}
\def\geq{\geqslant}

\def\1{{\bf 1}}

\def \<{\langle}
\def \>{\rangle}

\def \g{\mathfrak g}

\def \b{\beta}

\newenvironment{pf}{\begin{proof}}{\end{proof}}

\title[Minimal singular vectors in types B and C]
{On minimal singular vectors of $V^\kappa(\mathfrak{g})$ for types $B$ and $C$}
\author{Cuipo Jiang}
\author{Tong Wu}
\address{Cuipo Jiang and Tong Wu, School of Mathematical Sciences,
Shanghai Jiao Tong University, Shanghai, China}
\email[Jiang]{cpjiang@sjtu.edu.cn}
\email[Wu]{wt\_lovetju@sjtu.edu.cn}
\date{}

\begin{document}

\begin{abstract}
Let $\mathfrak{g}$ be  a finite-dimensional simple Lie algebra of type $B$ or $C$ over 
${\mathbb C}$. Given $\kappa\in{\mathbb C}$, let  $V^\kappa(\mathfrak{g})$  be the universal affine vertex algebra associated with $\mathfrak{g}$ at  level $\kappa$.  We  characterize the weights of minimal singular vectors of $V^\kappa(\mathfrak{g})$ with respect to the generalized Cartan matrix, when $V^\kappa(\mathfrak{g})$  is not simple and $\kappa$ is not critical. Furthermore, we determine  the  corresponding longest Weyl elements in the  Kashiwara-Tanisaki character theorem of  the  singular vectors.
\end{abstract}

\maketitle
\tableofcontents
\newpage

\section{Introduction}\label{sec:introduction}

Let $\mathfrak{g}$  be a finite-dimensional simple Lie algebra over $\CC$ with the normalized bilinear form $(\cdot|\cdot)$ such that $(\alpha|\alpha)=2$ when $\alpha$ is  a long root. Let $\widehat{\mathfrak{g}}=\mathfrak{g}\otimes\CC[t,t^{-1}]\oplus \CC K$ be the associated affine Lie algebra \cite{kac1990infinite}. Let $h$ and $\mathbf{h}^\vee$ denote the Coxeter number and the dual Coxeter number of $\mathfrak{g}$, respectively. For $\kappa\in\CC$, let  $V^\kappa(\mathfrak{g})$ be the universal affine vertex algebra associated to $\mathfrak{g}$ and $\kappa$, and  $L_\kappa(\mathfrak{g})$ the simple quotient of it. Recall from \cite{GK07} that $V^\kappa(\mathfrak{g})$ is not simple if and only if $r^\vee(\kappa+\mathbf{h}^\vee)\in\mathbb{Q}_{\geq 0}\setminus\{\frac{1}{m}:m\in\ZZ_{\geq 1}\}$, where $r^\vee$ is the lacing number of $\mathfrak{g}$.  It is known that when $\kappa$ is an admissible rational number, that is, 
	\begin{align*}
	\kappa+\mathbf{h}^\vee=\frac{p}{q},\quad p,q\in \ZZ_{\geq 1}, \ (p,q)=1,\
	p\geq \begin{cases}\mathbf{h}^\vee&\text{if }(r^{\vee},q)=1,\\
		h&\text{if }(r^{\vee},q)\ne 1,
	\end{cases}
\end{align*}
the maximal ideal of $V^\kappa(\mathfrak{g})$ is generated by one singular vector  \cite{kac1990infinite, FZ92,   KW89, LL04,Pe07,Pe08}, etc. When $\kappa=-\mathbf{h}^\vee+\frac{p}{q} $ is not admissible, there has been progress in determining   maximal ideals of $V^\kappa(\mathfrak{g})$ for some special cases. Specifically,   the maximal ideal of $V^\kappa(\mathfrak{g})$ is generated by two singular vectors when $\mathfrak{g}=\mathfrak{sl}_3$ and $\kappa+3=\frac{2}{2m+1}$ \cite{AP08, AKMPP20, JS25} .  While $\mathfrak{g}=\mathfrak{sl}_n$, $\kappa+n=\frac{n-1}{q}$, $n\geq 4$, and $q\geq 1$, the maximal ideal of $V^\kappa(\mathfrak{g})$ is generated by one singular vector \cite{AP08, AM17,APV23,JS25}.  For the case that $\mathfrak{g}$ is of $G_2$ type  and $\kappa=-2$, the  maximal ideal of $V^\kappa(\mathfrak{g})$ is generated by one singular vector \cite{ADFLM25}.  In general,  results on  singular vectors and maximal ideals remain  limited.  

Let $\mathfrak{h}$ and $\widehat{\mathfrak{h}}$ be the Cartan subalgebras of $\mathfrak{g}$ and $\widehat{\mathfrak{g}}$, respectively. Notice that  if $V^\kappa(\mathfrak{g})$ is not simple, it has at least one singular vector with minimal conformal weights which we will  call minimal singular vectors for simplicity.  When $\mathfrak{g}$ is of simply-laced type,  the weights of all minimal singular vectors  with respect to the Cartan subalgebra $\widehat{\mathfrak{h}}$ 
were determined using Gorelik-Kac determinant formula \cite{JS26}.  We treat the cases in which $\mathfrak{g}$ is of type $B$ or $C$.  Since $B$ and $C$ are not symmetric type and the associated affine Weyl groups  are more complicated, we need to refine the analysis process and consider more special cases. By Theorems  \ref{Bconclusion} and \ref{specialcaseB}, we give the  weights of all minimal singular vectors of $V^\kappa(\mathfrak{g})$. While by Theorems  \ref{C_conclusion} and \ref{specialcaseC} we give  the results for type  $C$ case.

Let $W$ be the Weyl group of $\mathfrak{g}$, and 
$\widehat{W}$  the affine Weyl group of $\widehat{\mathfrak{g}}$. For $\lambda\in\widehat{\mathfrak{h}}^*$, denote $\widehat\Delta(\lambda)=\{\alpha\in\widehat\Delta^{\mathrm{re}}|\langle\lambda+\widehat\rho,\alpha^\vee\rangle\in{\ZZ}\}$.  Let $\widehat{W}(\lambda)$ be the associated affine Weyl group generated by reflections $s_{\alpha}$, $\alpha\in \widehat\Delta(\lambda)$, and $\widehat{W}_0(\lambda)=\{w\in \widehat{W}(\lambda)|w\circ\lambda=\lambda\}$ the subgroup of 
$\widehat{W}(\lambda)$, where $w\circ \lambda=w(\lambda+\widehat\rho)-\widehat\rho$.  By the Kashiwara-Tanisaki character theorem \cite{KT00}, for given 
$\kappa\Lambda_0$ and the weight   $\lambda_{\nu}=\kappa\Lambda_0-qD_p\delta+\nu$ of a minimal singular vector in $V^{\kappa}(\mathfrak{g})$, there exist the unique corresponding dot-dominant $\Lambda\in{\widehat{\mathfrak h}}^*$ and longest  elements $\widetilde{y}$ in $\widetilde{y}\widehat{W}_0(\Lambda)$ and $\widetilde{z}=\widetilde{z}_0\widetilde{y}$ in  $\widetilde{z}\widehat{W}_0(\Lambda)$, respectively, 
such that  $$\kappa\Lambda_0=\widetilde{y}\circ\Lambda,$$
and 
$$\lambda_{\nu}=\widetilde{z}\circ \Lambda=\widetilde{z}_0\widetilde{y}\circ \Lambda.$$
\cite{JS26} characterizes $\widetilde{y}$ and $\widetilde{z}$ explicitly when $\mathfrak{g}$ is of ADE types. This paper deals with the cases that $\mathfrak{g}$ is of types $B$ and $C$. By Theorem \ref{thm:A-cosets} we give our results on $\widetilde{y}$ and $\widetilde{z}$.

The  rest of the paper is organized as follows.  In Section 2, we give some preliminaries.  Minimal conformal weights of maximal ideals of $V^\kappa(\mathfrak{g})$ are studied in Section 3.  We give results on weights of minimal singular vectors of  $V^\kappa(\mathfrak{g})$  when $\mathfrak{g}$ is of type $B$ in Section 4. Section 5 is focused on the case that $\mathfrak{g}$ is of type $C$. In Section 6, we characterize $\widetilde{y}$ and $\widetilde{z}$ when $\mathfrak{g}$ is of type $B$ or $C$.

\section{Preliminaries}\label{sec:preliminaries}
\subsection{Affine Lie algebras and affine vertex algebras}
Let $\mathfrak g$ be a finite-dimensional simple Lie algebra over $\CC$.
Fix a triangular decomposition
\(\mathfrak g=\mathfrak n_-\oplus\mathfrak h\oplus\mathfrak n_+\),
and let $\Delta=\Delta^+\sqcup\Delta^-$ be the corresponding root system.
Write \(Q=\sum_{\alpha\in\Delta}\ZZ\alpha\) for the root lattice. Let $W$ be the Weyl
group, $\ell$ its length function, and $\rho$ half of the sum of the positive
roots. We use the dot action \(w\circ\lambda=w(\lambda+\rho)-\rho\).
Normalize the invariant form by
\((\cdot|\cdot)=(2\mathbf{h}^\vee)^{-1}\times\)Killing form,
so every long root has squared length $2$. We use this form to identify
$\mathfrak h$ with $\mathfrak h^*$, and write
$\|\lambda\|^2=(\lambda|\lambda)$. The dual Coxeter number and the
lacing number are denoted by $\mathbf{h}^\vee$ and $r^\vee$, respectively.

The affine Lie algebra is
\[
\widehat{\mathfrak g}
=\mathfrak g\otimes\CC[t,t^{-1}]\oplus\CC K,
\]
with bracket
\[
[x\otimes t^m,y\otimes t^n]
=[x,y]\otimes t^{m+n}
+m\delta_{m,-n}(x|y)K.
\]
We write $x(m)=x\otimes t^m$ and regard affine roots and weights as formal elements of
\(\mathfrak h^*\oplus\CC\Lambda_0\oplus\CC\delta\).
Here $\Lambda_0|_{\mathfrak h}=0$ and $\Lambda_0(K)=1$. The
normalized form extends by
\[
(\delta|\delta)=(\Lambda_0|\Lambda_0)=0,
\qquad
(\delta|\Lambda_0)=1,
\qquad
(\delta|\mathfrak h^*)=(\Lambda_0|\mathfrak h^*)=0.
\]
The real and imaginary affine roots are
\[
\widehat\Delta^{\mathrm{re}}
=\{\alpha+m\delta:\alpha\in\Delta,\ m\in\ZZ\},
\qquad
\widehat\Delta^{\mathrm{im}}
=\{m\delta:m\in\ZZ\setminus\{0\}\}.
\]
Set
\(
\widehat\Delta
=\widehat\Delta^{\mathrm{re}}\sqcup\widehat\Delta^{\mathrm{im}}
\). The corresponding root spaces are denoted by
\[
\widehat{\mathfrak g}_{\alpha+m\delta}
:=\mathfrak g_\alpha\otimes t^m
\quad(\alpha\in\Delta,\ m\in\ZZ),
\qquad
\widehat{\mathfrak g}_{m\delta}
:=\mathfrak h\otimes t^m
\quad(m\ne0).
\]
We take
\[
\widehat\Delta^+
=\Delta^+
\cup\{\alpha+m\delta:\alpha\in\Delta,\ m>0\}
\cup\{m\delta:m>0\}.
\]
If
$\beta=\alpha+m\delta\in\widehat\Delta^{\mathrm{re}}$, then
\[
\beta^\vee
=\alpha^\vee+\frac{2m}{(\alpha|\alpha)}K,
\qquad
\alpha^\vee=\frac{2\alpha}{(\alpha|\alpha)}.
\]
We also set
\[
\widehat Q=Q\oplus\ZZ\delta,
\qquad
\widehat Q_+
=\sum_{\beta\in\widehat\Delta^+}\ZZ_{\geq0}\beta.
\]
Put
\[
\widehat{\mathfrak g}_+
=\mathfrak g\otimes t\CC[t],
\qquad
\widehat{\mathfrak n}_+
=\mathfrak n_+\oplus\widehat{\mathfrak g}_+.
\]
Let $\mathfrak g[t]=\mathfrak g\otimes\CC[t]$. For $\kappa\in\CC$, let
$\CC_\kappa$ be the one-dimensional
$\mathfrak g[t]\oplus\CC K$-module on which $\mathfrak g[t]$ acts
trivially and $K$ acts by $\kappa$. The vacuum module is
\[
V^\kappa(\mathfrak g)
=U(\widehat{\mathfrak g})
\otimes_{U(\mathfrak g[t]\oplus\CC K)}\CC_\kappa.
\]
Its generator $\mathbf1$ is the vacuum vector of affine highest
weight $\kappa\Lambda_0$. Let
\[
\widehat Q_+^{\mathrm{vac}}
:=\sum_{\gamma\in\widehat\Delta^+\setminus\Delta}
\ZZ_{\geq0}\gamma
\subset\widehat Q_+.
\]
For \(\nu\in\widehat Q_+^{\mathrm{vac}}\), we write
\(V^\kappa(\mathfrak g)_\nu\) for the weight space of affine weight
\(\kappa\Lambda_0-\nu\), and the same convention applies to submodules and
subquotients.
The vacuum module carries the standard grading
\[
V^\kappa(\mathfrak g)
=\bigoplus_{c\geq0}V^\kappa(\mathfrak g)[c],
\qquad
\deg\mathbf1=0,
\qquad
\deg x(-m)=m\quad(m\geq1).
\]
 It is known that  $V^\kappa(\kg)$ carries a structure of vertex algebra, called the universal affine vertex algebra associated to $\mathfrak{g}$ at level $\kappa$\cite{FZ92,LL04}.
 
Choose dual bases \(\{u_i\}\) and \(\{u^i\}\) of \(\mathfrak g\)
with respect to \((\cdot|\cdot)\). For \(\kappa\ne-\mathbf{h}^\vee\), the Sugawara operator \(L(0)\) on the
vacuum module is given by
\[
L(0)
=\frac{1}{2(\kappa+\mathbf{h}^\vee)}
\left(
\sum_i u_i(0)u^i(0)
+2\sum_{m\geq1}\sum_i u_i(-m)u^i(m)
\right).
\]
The relations \([L(0),x(m)]=-m x(m)\) and \(L(0)\mathbf1=0\)
show that \(V^\kappa(\mathfrak g)[c]\) is the eigenspace of
\(L(0)\) with eigenvalue \(c\). Thus the standard grading
coincides with the Sugawara conformal weight grading.
Only finitely many terms in the displayed sum act nontrivially
on each vector. Hence every affine submodule is preserved by
\(L(0)\) and is graded.
The unique maximal proper submodule is denoted by \(N_\kappa(\mathfrak g)\).
For a graded subquotient \(M\), we write \(M[c]\) for its
subspace of conformal weight \(c\).
A nonzero affine weight vector \(v\in N_\kappa(\mathfrak g)\)
is called a singular vector if \(\widehat{\mathfrak n}_+v=0\).
The reducibility criterion takes the following form \cite{GK07}.

\begin{theorem}
\label{nonsimple}
The vacuum module \(V^\kappa(\mathfrak g)\) is reducible
if and only if
\[
r^\vee(\kappa+\mathbf{h}^\vee)
\in\mathbb Q_{\geq0}
\setminus\left\{\frac1m:m\in\ZZ_{\geq1}\right\}.
\]
\end{theorem}

\subsection{Shapovalov form and Jantzen filtration}
Let $\sigma_0$ be the Chevalley anti-involution of $\mathfrak g$, so that
$\sigma_0|_{\mathfrak h}=\mathrm{id}$ and
$\sigma_0(\mathfrak g_\alpha)=\mathfrak g_{-\alpha}$. Extend it to the
affine algebra by
\(\sigma(x(m))=\sigma_0(x)(-m)\) and \(\sigma(K)=K\).
There is a unique normalized symmetric contravariant form\cite{shapovalov1972bilinear,KK79}
\[
S_\kappa:V^\kappa(\mathfrak g)\times V^\kappa(\mathfrak g)\longrightarrow\CC
\]
satisfying
\[
S_\kappa(\mathbf1,\mathbf1)=1,
\qquad
S_\kappa(uv_1,v_2)=S_\kappa(v_1,\sigma(u)v_2).
\]
The radical of \(S_\kappa\) is \(N_\kappa(\mathfrak g)\). Set
\(\widehat\rho=\rho+\mathbf{h}^\vee\Lambda_0\). Define
\[
\widehat R
=\prod_{\alpha\in\widehat\Delta^+}
(1-e^{-\alpha})^{\dim\widehat{\mathfrak g}_\alpha}.
\]
For $\nu\in\widehat Q_+$, let $S_\nu(\kappa)$ be the restriction of
$S_\kappa$ to $V^\kappa(\mathfrak g)_\nu$. Specializing the vacuum determinant formula of
\cite{GK07} to a finite-dimensional
simple Lie algebra, and using the normalization fixed above, gives
\[
\det S_\nu(\kappa)
=\prod_{r=1}^{\infty}
\prod_{\gamma\in\widehat\Delta^+\setminus\Delta}
\phi_{r,\gamma}(\kappa)^{
d_{r,\gamma}(\nu)\dim\widehat{\mathfrak g}_\gamma},
\]
where
\[
\phi_{r,\gamma}(\kappa)
=(\kappa\Lambda_0+\widehat\rho|\gamma)
-\frac r2(\gamma|\gamma),
\]
and the integers $d_{r,\gamma}(\nu)$ are determined by
\[
\sum_{\nu\in\widehat Q_+}d_{r,\gamma}(\nu)e^{-\nu}
=\widehat R^{-1}
\sum_{w\in W}(-1)^{\ell(w)}e^{w\circ(-r\gamma)}.
\]
For fixed \(\nu\), only finitely many exponents
\(d_{r,\gamma}(\nu)\dim\widehat{\mathfrak g}_\gamma\) are nonzero.
Hence there is a nonzero constant $C_\nu$ and there are nonnegative integers $m_b(\nu)$,
only finitely many of which are nonzero, such that
\[
\det S_\nu(\kappa)
=C_\nu
\prod_{\substack{b\in\CC\\m_b(\nu)>0}}
(\kappa+\mathbf{h}^\vee-b)^{m_b(\nu)}.
\]
Thus
\(N_\kappa(\mathfrak g)_\nu\ne0\) if and only if
\(m_{\kappa+\mathbf{h}^\vee}(\nu)>0\).
Let $A=\CC[t]_{(t)}$ be the localization of ${\mathbb C}[t]$ by the maximal ideal  generated by $t$, and deform the level to
$\kappa+t$. Denote by $V_A^{\kappa+t}(\mathfrak g)$ 
the  vacuum module with the corresponding $A$-form, 
and by
$S_{\kappa+t}$ its contravariant form. For $i\geq0$, set
\[
V_A^i
=\{v\in V_A^{\kappa+t}(\mathfrak g):
S_{\kappa+t}(v,V_A^{\kappa+t}(\mathfrak g))
\subset t^iA\}.
\]
After identifying
$V_A^{\kappa+t}(\mathfrak g)/t
V_A^{\kappa+t}(\mathfrak g)$ with $V^\kappa(\mathfrak g)$, define
\[
\mathcal F^i(V^\kappa(\mathfrak g))
=(V_A^i+t V_A^{\kappa+t}(\mathfrak g))/
t V_A^{\kappa+t}(\mathfrak g).
\]
This is the Jantzen filtration. It is decreasing  and
\[
\mathcal F^0(V^\kappa(\mathfrak g))=V^\kappa(\mathfrak g),
\qquad
\mathcal F^1(V^\kappa(\mathfrak g))=N_\kappa(\mathfrak g).
\]
For $\nu\in\widehat Q_+$, set
\[
\mathcal F^i(V^\kappa(\mathfrak g))_\nu
=\mathcal F^i(V^\kappa(\mathfrak g))
\cap V^\kappa(\mathfrak g)_\nu.
\]
On every weight space, we have the Jantzen sum formula
\cite{GK07,Jan77}
\begin{equation}\label{sumformula}
\sum_{i=1}^{\infty}
\dim\mathcal F^i(V^\kappa(\mathfrak g))_\nu
=
\mathrm{ord}_{t=0}\det S_\nu(\kappa+t)
=m_{\kappa+\mathbf{h}^\vee}(\nu).
\end{equation}

\section{\texorpdfstring{Minimal conformal weights of
$N_{\kappa}(\mathfrak g)$}{Minimal conformal weight of the maximal submodule}}
\label{sec:minimal-weight}
For a subquotient \(M\) of \(V^\kappa(\mathfrak g)\), define its
normalized character by
\(\operatorname{ch}M
=\sum_{\nu\in\widehat Q_+}(\dim M_\nu)e^{-\nu}\),
where \(M_\nu\) is the space of affine weight
\(\kappa\Lambda_0-\nu\).
Summing \eqref{sumformula} over \(\nu\in\widehat Q_+\) gives
\[
\sum_{i=1}^{\infty}
\operatorname{ch}\mathcal F^i(V^\kappa(\mathfrak g))
=\widehat R^{-1}M_\kappa,
\]
where
\[
M_\kappa
:=\widehat R\sum_{\nu\in\widehat Q_+}
m_{\kappa+\mathbf{h}^\vee}(\nu)e^{-\nu}.
\]
Since the filtration is decreasing and
\(\mathcal F^1(V^\kappa(\mathfrak g))=N_\kappa(\mathfrak g)\),
this identity shows that \(V^\kappa(\mathfrak g)\) is reducible
if and only if \(M_\kappa\ne0\).
Assume that \(\kappa\ne-\mathbf{h}^\vee\) and
\(V^\kappa(\mathfrak g)\) is reducible.
By Theorem~\ref{nonsimple}, we may write
\(r^\vee(\kappa+\mathbf{h}^\vee)=p/q\), where
\(p\in\ZZ_{\ge2}\), \(q\in\ZZ_{\ge1}\), and \((p,q)=1\).

Set
\[
E^\rho(\lambda)
=\sum_{w\in W}(-1)^{\ell(w)}e^{w\circ\lambda}.
\]
We also write
\[
E(\mu)
:=\sum_{w\in W}(-1)^{\ell(w)}e^{w(\mu)-\rho}
=E^\rho(\mu-\rho).
\]
The determinant formula \cite{GK07} gives
\[
M_\kappa
=\sum_{\substack{
r\in\ZZ_{\ge1},\
\gamma\in\widehat\Delta^+\setminus\Delta\\
\phi_{r,\gamma}(\kappa)=0}}
\dim\widehat{\mathfrak g}_\gamma\,E^\rho(-r\gamma).
\]
Since \(\phi_{r,m\delta}(\kappa)=m(\kappa+\mathbf{h}^\vee)\ne0\),
only real roots contribute to this sum.
Write \(\gamma=m\delta+\alpha\), where
\(m\in\ZZ_{\ge1}\) and \(\alpha\in\Delta\).
Then \(\dim\widehat{\mathfrak g}_{m\delta+\alpha}=1\), and
\[
\phi_{r,m\delta+\alpha}(\kappa)
=m(\kappa+\mathbf{h}^\vee)+(\rho|\alpha)
-\frac r2\|\alpha\|^2.
\]
Thus \(\phi_{r,m\delta+\alpha}(\kappa)=0\) is equivalent to
\[
mp
=q\,r^\vee\left(\frac r2\|\alpha\|^2-(\rho|\alpha)\right).
\]
Together with \((p,q)=1\), this implies \(m=qs\) for some
\(s\in\ZZ_{\ge1}\).
The equation then becomes
\[
r^\vee\left(\frac r2\|\alpha\|^2-(\rho|\alpha)\right)=ps.
\]
For \(r,s\in\ZZ_{\ge1}\), set
\begin{equation}\label{3ce1}
\Delta_{r,s}^p
:=\left\{\alpha\in\Delta:
r^\vee\left(\frac r2\|\alpha\|^2-(\rho|\alpha)\right)=ps
\right\},
\end{equation}
and for \(D\in\ZZ_{\ge1}\), set
\[
M_{p,D}
:=\sum_{\substack{r,s\in\ZZ_{\ge1}\\rs=D}}
\sum_{\alpha\in\Delta_{r,s}^p}E^\rho(-r\alpha).
\]
It follows that
\[
M_\kappa
=\sum_{D\in\ZZ_{\ge1}}e^{-qD\delta}M_{p,D}.
\]
Let \(D_p\) be the least positive integer \(D\) for which
\(M_{p,D}\ne0\).
The Jantzen sum formula then shows that the minimal conformal
weight of \(N_\kappa(\mathfrak g)\) is \(qD_p\).

For a dominant integral weight \(\lambda\) of $\g$, let
\(L_{\mathfrak g}(\lambda)\) denote the simple finite-dimensional
\(\mathfrak g\)-module of highest weight \(\lambda\).

Put \(d=qD_p\) and
\(R=\prod_{\alpha\in\Delta^+}(1-e^{-\alpha})\).
For a homogeneous component \(M[c]\), write
\(\operatorname{ch}M[c]\) for its character as a \(\mathfrak g\)-module.
Taking the coefficient of \(e^{-d\delta}\) in the Jantzen sum
formula gives
\begin{equation}\label{eq:first-coefficient-positivity}
R^{-1}M_{p,D_p}
=\sum_{j\ge1}
\operatorname{ch}\mathcal F^j(V^\kappa(\mathfrak g))[d].
\end{equation}
By complete reducibility and the Weyl character formula
\(R\operatorname{ch}L_{\mathfrak g}(\nu)=E^\rho(\nu)\), the element
\(M_{p,D_p}\) is a nonnegative integral linear combination of
\(E^\rho(\nu)\), where \(\nu\) ranges over dominant integral weights.

\begin{prop}[{\cite[Theorem~3.4]{JS26}}]
\label{prop:first-jantzen-layer}
Suppose that
\begin{align}
M_{p,D_p}=\sum_{i=1}^m a_iE^\rho(\lambda_i),
\end{align}
where \(a_i\in\ZZ_{>0}\) and
\(\lambda_1,\ldots,\lambda_m\) are pairwise distinct dominant
integral weights.
Then \(N_\kappa(\mathfrak g)[c]=0\) for \(0\le c<qD_p\), and
\(N_\kappa(\mathfrak g)[qD_p]
\cong\displaystyle\bigoplus_{i=1}^m
L_{\mathfrak g}(\lambda_i)^{\oplus a_i}\)
as a \(\mathfrak g\)-module.
At conformal weight \(qD_p\), the space of singular vectors of
weight \(\kappa\Lambda_0-qD_p\delta+\lambda_i\) has dimension
\(a_i\). 
\end{prop}

\begin{lem}\label{general E=0 criteria}
For \(w\in W\), one has
\[
E^\rho(w\circ\lambda)=(-1)^{\ell(w)}E^\rho(\lambda).
\]
If \(E^\rho(\lambda_1)\) and \(E^\rho(\lambda_2)\) are nonzero,
then \(E^\rho(\lambda_1)=\pm E^\rho(\lambda_2)\) if and only if
\(\lambda_1=w\circ\lambda_2\) for some \(w\in W\).
Moreover, \(E^\rho(\lambda)=0\) if and only if
\(\lambda+\rho\) is fixed by a reflection in \(W\).
\end{lem}

\begin{pf}
The first identity follows by reindexing.
If \(\lambda+\rho\) is regular, the exponents
\(w(\lambda+\rho)-\rho\) are pairwise distinct.
Otherwise, a reflection fixing \(\lambda+\rho\) pairs terms
with opposite signs.
The remaining assertion follows by comparing exponents.
\end{pf}

From now on, let \(\mathfrak g\) be of type \(X_n\), where
\(X=B\) with \(n\ge2\), or \(X=C\) with \(n\ge3\).
In both cases \(r^\vee=2\). Set
\[
\Delta_{\mathrm{short}}
=\{\alpha\in\Delta:\|\alpha\|^2=1\},
\qquad
\Delta_{\mathrm{long}}
=\{\alpha\in\Delta:\|\alpha\|^2=2\}.
\]
Write \(\lambda=(a_1,\ldots,a_n)\) when
\(\lambda=\sum_i a_i e_i\).
In type \(B_n\), one has
\(\rho=(n-\frac12,n-\frac32,\dots,\frac12)\),
with short roots \(\pm e_i\) and long roots
\(\pm(e_i\pm e_j)\), \(i<j\).
In type \(C_n\), one has \(\rho=(n,n-1,\dots,1)\),
with short roots \(\pm(e_i\pm e_j)\), \(i<j\), and long roots
\(\pm2e_i\).
The Weyl group permutes the coordinates and changes their
signs independently. Thus the dominant representative of
\(\rho-a\alpha\) is obtained by arranging the absolute values
of its coordinates in nonincreasing order.

\begin{lem}\label{lem:signed-permutation-alternants}
For \(\mu=(a_1,\ldots,a_n)\), one has \(E(\mu)=0\) if and only if
\(a_i=0\) for some \(i\), or \(a_i=\pm a_j\) for some \(i\ne j\).
\end{lem}

\begin{pf}
Apply Lemma~\ref{general E=0 criteria} with \(\lambda=\mu-\rho\)
to the reflecting hyperplanes \(a_i=0\) and \(a_i=\pm a_j\).
\end{pf}

Let \(\mathcal S^X_{p,n}\) denote the set of dominant integral
weights \(\nu\) for which the coefficient of \(E^\rho(\nu)\)
in \(M_{p,D_p}\) is nonzero.
For \(a\in\ZZ_{\ge1}\) and \(\alpha\in\Delta\) with
\(E^\rho(-a\alpha)\ne0\), let \(\mu\) be the dominant
representative of \(\rho-a\alpha\) and put \(\nu=\mu-\rho\).
Then \(\nu\) is dominant integral, and
\(E^\rho(-a\alpha)=\pm E(\mu)=\pm E^\rho(\nu)\).

Recall  $\Delta_{r,s}^p$ from (\ref{3ce1}), then we have
\(\Delta_{r,s}^p
=\Delta_{r,s}^{p,\mathrm{short}}
\sqcup\Delta_{r,s}^{p,\mathrm{long}}\), where
\[
\begin{aligned}
\Delta_{r,s}^{p,\mathrm{short}}
&=\{\alpha\in\Delta_{\mathrm{short}}:
r-2(\rho|\alpha)=ps\},\\
\Delta_{r,s}^{p,\mathrm{long}}
&=\{\alpha\in\Delta_{\mathrm{long}}:
2(r-(\rho|\alpha))=ps\}.
\end{aligned}
\]
Set
\begin{equation}\label{ce1}
P_{p,D}
:=\{(a,b)\in\ZZ_{\ge1}^2:ab=pD,\ p\mid b\}.
\end{equation}
For \(A\subset\ZZ_{\ge1}^2\), define
\[
A':=\{(b,a):(a,b)\in A\},
\qquad
t(A):=\left\{\left(2a,\frac b2\right):
(a,b)\in A,\ 2\mid b\right\}.
\]
For pairs, write \((a,b)'=(b,a)\), and  when \(b\) is even,
\(t(a,b)=(2a,b/2)\).
For fixed \(p,D\), write
\[
P=P_{p,D},
\qquad
U=P\mathbin{\triangle}P'
=(P\setminus P')\cup(P'\setminus P).
\]
Substituting \(a=r\) and \(b=ps\) gives
\[
M_{p,D}=\sum_{(a,b)\in P_{p,D}}E_{a,b},
\]
where
\[
\Omega_{a,b}
:=
\{\alpha\in\Delta_{\mathrm{short}}:
2(\rho|\alpha)=a-b\}
\sqcup
\{\alpha\in\Delta_{\mathrm{long}}:
2(\rho|\alpha)=2a-b\},
\]
and
\[
E_{a,b}
:=\sum_{\alpha\in\Omega_{a,b}}E^\rho(-a\alpha).
\]

Let \(\alpha\in\Omega_{a,b}\cap\Delta^-\).
If \(\alpha\) is short, then \(-\alpha\in\Omega_{b,a}\).
If \(\alpha\) is long, then
\(b=2a-2(\rho|\alpha)\) is even and
\(-\alpha\in\Omega_{b/2,2a}\).
The defining equations give
\[
\begin{aligned}
s_\alpha(\rho-a\alpha)&=\rho+b\alpha
&&(\alpha\in\Delta_{\mathrm{short}}),\\
s_\alpha(\rho-a\alpha)&=\rho+\frac b2\alpha
&&(\alpha\in\Delta_{\mathrm{long}}).
\end{aligned}
\]
Lemma~\ref{general E=0 criteria} therefore gives
\[
\begin{aligned}
E^\rho(b\alpha)&=-E^\rho(-a\alpha)
&&(\alpha\in\Delta_{\mathrm{short}}),\\
E^\rho\left(\frac b2\alpha\right)&=-E^\rho(-a\alpha)
&&(\alpha\in\Delta_{\mathrm{long}}).
\end{aligned}
\]

The coefficients \(M_{p,D}\), and hence \(D_p\)  are independent
of \(q\).
The admissible ranges \cite{KW89} and the corresponding minimal
conformal weights of singular vectors
\cite[Corollary~2.1]{CL22} are
\[
\begin{array}{c|c|c}
\text{type}
&p\text{ even: principal}
&p\text{ odd: coprincipal}\\ \hline
B_n&
p\ge4n-2,\quad q\left(\dfrac p2-2n+2\right)&
p\ge2n,\quad q(p-2n+1)\\[2mm]
C_n&
p\ge2n+2,\quad q\left(\dfrac p2-n\right)&
p\ge2n,\quad q(p-2n+1).
\end{array}
\]
These are the first two cases in
Theorems~\ref{Bconclusion} and~\ref{C_conclusion}.
The remaining  are nonadmissible cases.

\begin{lem}
\label{lem:signed-permutation-orbit-recovery}
Write \(\rho=(r_1,\ldots,r_n)\) with
\(r_1>\cdots>r_n>0\).
\begin{enumerate}
\item
If a dominant representative \(\mu\) is obtained from \(\rho\)
by deleting a coordinate \(v\) and adjoining \(u>r_1\),
then \(\mu\) uniquely determines \(u,v\).
\item
If a dominant representative \(\mu\) is obtained from \(\rho\)
by deleting \(r_i,r_j\), with \(1\le i<j\le n\), and adjoining
\[
\{r_i+a,r_j+a\}
\quad\text{or}\quad
\{r_i+a,a-r_j\},
\]
where \(a\in\ZZ_{\ge1}\) and both new coordinates exceed \(r_1\),
then \(\mu\) uniquely determines \(i,j,a\) and which pair occurs.
\item
No representative in part~(1) equals a representative in part~(2).
\end{enumerate}
For distinct representatives \(\mu\) in parts~(1) and~(2),
the corresponding \(E(\mu)\) are nonzero and linearly independent.
\end{lem}

\begin{pf}
In part~(1), \(u\) is the largest coordinate of \(\mu\),
and \(v\) is the unique coordinate of \(\rho\) absent from \(\mu\).
In part~(2), the missing coordinates determine \(i,j\),
and subtracting \(r_i\) from the largest coordinate gives \(a\).
The two largest coordinates differ by \(r_i-r_j\) or
\(r_i+r_j\), respectively, so the two cases are distinct.
Part~(3) follows by counting the coordinates greater than \(r_1\).
All these representatives are regular.
For distinct dominant representatives, the corresponding
\(E(\mu)\) have disjoint sets of exponents and are therefore
linearly independent.
\end{pf}

The next two sections give formulas for \(M_{p,D}\).
The values of \(D_p\) and the coefficients \(M_{p,D_p}\)
are determined in Sections~\ref{sec:type-B}
and~\ref{sec:type-C}.

\section{\texorpdfstring{Type $B_n$}{Type Bn}}\label{sec:type-B}
Fix \(n\geq2\). Let
$\mathfrak g$ be of type $B_n$. We use the standard coordinate realization of the root system
of type \(B_n\) \cite{Humphreys,tauvel2005lie} as
\[
\alpha_i=e_i-e_{i+1}\quad(1\leq i<n),\qquad \alpha_n=e_n,
\]
where $e_1,\dots,e_n$ is the standard orthonormal basis. Thus
\[
\Delta_{\mathrm{short}}=\{\pm e_i:1\leq i\leq n\},\qquad
\Delta_{\mathrm{long}}=\{\pm(e_i\pm e_j):1\leq i<j\leq n\},
\]
and for fixed $p$ and $D$ we write \(P=P_{p,D}\), so that
\(M_{p,D}=\sum_{(a,b)\in P}E_{a,b}\).

For \(a\in\ZZ_{\ge1}\) and \(\alpha\in\Delta\) with
\(E^\rho(-a\alpha)\ne0\), let \(\mu\) be the dominant representative
of \(\rho-a\alpha\).
We say that \(E^\rho(-a\alpha)\) is of type~\(\mathrm I\)
if \(\mu\) has exactly one coordinate greater than
\(n-\tfrac12\), and of type~\(\mathrm{II}\) if \(\mu\)
has exactly two such coordinates.
\begin{lem}\label{nonzero E criteria for B complex}
Let \(a\in\ZZ_{\ge1}\) and \(\alpha\in\Delta^-\).
Then \(E^\rho(-a\alpha)\ne0\) precisely in the following cases.
Here \(1\le i\le n\) for \(\alpha=-e_i\), and
\(1\le i<j\le n\) for roots involving two indices.
\begin{enumerate}[(i)]
	\item Type $\mathrm{I}$:
	\begin{enumerate}[(a)]
\item \(\alpha=-e_i\) and \(a\ge i\).
\item \(\alpha=-(e_i+e_j)\) and \(a=j-i\ge i\).
\item \(\alpha=-(e_i-e_j)\) and \(a=2n-i-j+1\ge i\).
\item \(\alpha=-(e_i-e_j)\) and \(a=2n-2j+1\ge i\).
\end{enumerate}
The dominant representative in cases~(a) and~(d) is
\[
(n-i+\tfrac12+a,n-\tfrac12,\dots,\widehat{n-i+\tfrac12},\dots,\tfrac12),
\]
and in cases~(b) and~(c) it is
\[
(n-i+\tfrac12+a,n-\tfrac12,\dots,\widehat{n-j+\tfrac12},\dots,\tfrac12).
\]

	\item Type $\mathrm{II}$:
	\begin{enumerate}[(a)]
		\item $\alpha=-(e_i+e_j)$ and $a\geq j$, with dominant representative
		\[
		(n-i+\frac12+a,n-j+\frac12+a,n-\frac12,\dots,\widehat{n-i+\frac12},\dots,\widehat{n-j+\frac12},\dots,\frac12).
		\]
		\item $\alpha=-(e_i-e_j)$ and $a\geq 2n-j+1$, with dominant representative
		\[
		(n-i+\frac12+a,a-(n-j+\frac12),n-\frac12,\dots,\widehat{n-i+\frac12},\dots,\widehat{n-j+\frac12},\dots,\frac12).
		\]
	\end{enumerate}
\end{enumerate}

\end{lem}
\begin{pf}
The negative roots are \(-e_i\), \(-(e_i+e_j)\), and
\(-(e_i-e_j)\). For each negative root \(\alpha\), the coordinates of
\(\rho-a\alpha\) are half-integers. Then by
Lemma~\ref{lem:signed-permutation-alternants} regularity amounts to
distinct absolute values.
For \(\alpha=-e_i\), the changed coordinate avoids the unchanged
coordinates exactly when \(a\ge i\), giving case~(i)(a).

For \(\alpha=-(e_i+e_j)\), regularity first requires \(a\ge i\).
If \(a<j\), the second changed coordinate must equal the omitted
\(i\)-th coordinate. So \(a=j-i\). Otherwise both changed
coordinates exceed \(n-\tfrac12\). These are cases~(i)(b) and~(ii)(a).

For \(\alpha=-(e_i-e_j)\), again \(a\ge i\) is necessary.
The absolute value of the second changed coordinate must either equal
one of the omitted coordinates \(n-i+\tfrac12,n-j+\tfrac12\),
or exceed \(n-\tfrac12\). The first alternatives give
\(a=2n-i-j+1\) and \(a=2n-2j+1\); while the last gives that
\(a\ge2n-j+1\). These are precisely the remaining three cases.
Sorting in each case gives the representatives in the statement.
\end{pf}
For \(\tau\in\{\mathrm{I},\mathrm{II}\}\), define
\[
\Omega_{a,b}^{\tau}
:=\{\alpha\in\Omega_{a,b}:
E^\rho(-a\alpha)\ne0
\text{ and }E^\rho(-a\alpha)\text{ is of type }\tau\},
\]
and put
\[
E_{a,b}^{\tau}
:=\sum_{\alpha\in\Omega_{a,b}^{\tau}}E^\rho(-a\alpha),
\qquad
E_P^{\tau}
:=\sum_{(a,b)\in P}E_{a,b}^{\tau}.
\]
Every nonzero term in \(E_{a,b}\) belongs to exactly one of these two types.
Consequently,
\begin{equation}\label{B-term-decomposition}
M_{p,D}
=E_P^{\mathrm{I}}
+E_P^{\mathrm{II}}.
\end{equation}
The nonzero terms described in Lemma~\ref{nonzero E criteria for B complex}
are indexed by the following subsets of \(\ZZ_{\geq1}^2\):
\begin{center}
\small
\(\displaystyle
\begin{aligned}
&S_{\mathrm{I},-}^{(1)}:=\{(a,b):1\leq b-a\leq 2n-1,\ b+a\geq 2n+1,\ 2\nmid(a+b)\},\\
&S_{\mathrm{I},-}^{(2)}:=\{(a,b):2\mid b,\ 2\nmid \tfrac b2,\ b\leq 4n-2,\ 4a+b\geq 4n+2,\ 2\leq b-4a\},\\
&S_{\mathrm{I},-}^{(3)}:=\{(a,b):2\mid b,\ 2\nmid \tfrac b2,\ 4a+b\geq 4n+2,\ 2\leq4a-b,\ b\leq 4n-2,\ 2a<b\},\\
&S_{\mathrm{I},-}^{(4)}:=\{(a,b):2\mid b,\ 2\nmid a,\ a+b\geq 2n+1,\ b-a\leq 2n-1,\ 2a<b\},\\
&S_{\mathrm{II},-}:=\{(a,b):2\mid b,\ 2\leq b-2a\leq 4n-4,\ b+2a\geq 4n+4\}.
\end{aligned}
\)
\end{center}
Let $S_{\mathrm{I},+}^{(1)}:=(S_{\mathrm{I},-}^{(1)})^{\prime}$,
$S_{\mathrm{II},+}:=(t(S_{\mathrm{II},-}))^{\prime}$, and
$S_{\mathrm{I},+}^{(j)}:=(t(S_{\mathrm{I},-}^{(j)}))^{\prime}$ for $j=2,3,4$. Put
$S_{\mathrm{I},\pm}=\bigcup_j S_{\mathrm{I},\pm}^{(j)}$. Further, define
\(S_{\mathrm{I}}=S_{\mathrm{I},+}\cup S_{\mathrm{I},-}\) and
\(S_{\mathrm{II}}=S_{\mathrm{II},+}\cup S_{\mathrm{II},-}\).
The sets \(S_{\mathrm I}\) and \(S_{\mathrm{II}}\) need not be disjoint. For example, in type \(B_3\),
\(E_{3,10}\) contains the type~\(\mathrm{I}\) term associated with
\(\alpha=-(e_1-e_3)\) and the type~\(\mathrm{II}\) term associated with
\(\alpha=-(e_2+e_3)\).

Put \(\Gamma=S_{\mathrm I,-}^{(1)}\).
For \((a,b)\in\Gamma\), put
\[
A_{a,b}
:=(-1)^{\,n-\frac12-\frac{b-a}{2}}
E\!\left(
\frac{a+b}{2},\,
n-\frac12,\,
\dots,\,
\widehat{\frac{b-a}{2}},\,
\dots,\,
\frac12
\right).
\]
Put \(j=n+\tfrac12-\tfrac{b-a}{2}\in\{1,\ldots,n\}\).
Then
\[
\nu=\left(\frac{a+b}{2}-n+\frac12\right)e_1
+\sum_{i=2}^{j}e_i,
\qquad
A_{a,b}=(-1)^{j-1}E^\rho(\nu),
\]
where the sum is empty when \(j=1\).

\begin{lem}\label{nonzero term set}
For \(a,b\in\ZZ_{\ge1}\) and
\(\tau\in\{\mathrm I,\mathrm{II}\}\), the sum \(E_{a,b}\)
contains a nonzero term of type \(\tau\) if and only if
\((a,b)\in S_\tau\).
\end{lem}

\begin{pf}
For type~\(\mathrm I\), substituting the four cases for negative
roots in Lemma~\ref{nonzero E criteria for B complex}
into the equations defining \(\Omega_{a,b}\) gives
\[
\begin{array}{c|c|c|c}
\text{case}&i&j&E^\rho(-a\alpha)\\ \hline
\text{(a)}&n-\frac{b-a}{2}+\frac12&-&A_{a,b}\\
\text{(b)}&n-\frac b4+\frac12&
n-\frac{b-4a}{4}+\frac12&-A_{t(a,b)}\\
\text{(c)}&n-\frac b4+\frac12&
n+\frac{b-4a}{4}+\frac12&A_{t(a,b)'}\\
\text{(d)}&n-\frac{b-a}{2}+\frac12&
n-\frac a2+\frac12&-A_{a,b}.
\end{array}
\]
Under these substitutions, the conditions on \(i,j,a\) in
cases~(a)--(d) are equivalent to those defining
\(S_{\mathrm{I},-}^{(1)},\ldots,S_{\mathrm{I},-}^{(4)}\),
respectively.
The signs follow from the coordinate permutations and sign
changes taking \(\rho-a\alpha\) to its dominant representative.

The term for \(-\alpha\) has the opposite sign and occurs at
\((b,a)\) in case~(a), or at \((b/2,2a)\) in cases~(b)--(d).
These pairs form the sets \(S_{\mathrm I,+}^{(j)}\).

For type~\(\mathrm{II}\) and \(\alpha\in\Delta^-\),
Lemma~\ref{nonzero E criteria for B complex} and the equations
defining \(\Omega_{a,b}\) give
\[
b=
\begin{cases}
2a+4n-2i-2j+2,
&\alpha=-(e_i+e_j),\quad a\ge j,\\
2a+2j-2i,
&\alpha=-(e_i-e_j),\quad a\ge2n-j+1.
\end{cases}
\]
In both cases, \(b\) is even,
\(2\le b-2a\le4n-4\), and \(b+2a\ge4n+4\).
Thus \((a,b)\in S_{\mathrm{II},-}\).

Conversely, let \((a,b)\in S_{\mathrm{II},-}\) and put
\(c=(b-2a)/2\).
Then \(c\) is an integer with \(1\le c\le2n-2\), and
\(2a+c\ge2n+2\).
If \(c\ge2\), take
\[
j=\left\lceil n+1-\frac c2\right\rceil,\qquad
i=2n+1-c-j,\qquad
\alpha=-(e_i+e_j).
\]
The bounds on \(c\) give \(1\le i<j\le n\), and
\(2a+c\ge2n+2\) gives \(a\ge j\).
Since \(2(\rho|\alpha)=-2c=2a-b\), we have
\(\alpha\in\Omega_{a,b}\).
Lemma~\ref{nonzero E criteria for B complex} therefore shows
that \(E^\rho(-a\alpha)\) is nonzero and of type~\(\mathrm{II}\).

If \(c=1\), take \(\alpha=-(e_{n-1}-e_n)\).
Here \(2a+1\ge2n+2\) implies \(a\ge n+1\), and
\(2(\rho|\alpha)=-2=2a-b\).
Thus \(\alpha\in\Omega_{a,b}\), and the same lemma gives a
nonzero term of type~\(\mathrm{II}\).

All roots contributing terms of type~\(\mathrm{II}\) are long roots.
Replacing \(\alpha\) by \(-\alpha\) sends \((a,b)\) to
\((b/2,2a)\) and reverses the sign of the corresponding term.
Hence the pairs for positive roots form
\(t(S_{\mathrm{II},-})'=S_{\mathrm{II},+}\).
\end{pf}

\begin{cor}\label{independent type I B}
The elements \(A_{a,b}\), with \((a,b)\in\Gamma\), are nonzero
and linearly independent.
\end{cor}

\begin{pf}
Put \(u=(a+b)/2\) and \(v=(b-a)/2\). Since \((a,b)\in\Gamma\),
\[
u\ge n+\frac12,\qquad
v\in\left\{\frac12,\frac32,\ldots,n-\frac12\right\}.
\]
The coordinates of the dominant representative are \(u\) and the
coordinates of \(\rho\) other than \(v\). Since \(a=u-v\) and \(b=u+v\), the assertion
follows from Lemma~\ref{lem:signed-permutation-orbit-recovery}(1).
\end{pf}

\begin{cor}\label{independent type II B}
Let \(a\in\ZZ_{\ge1}\) and \(\alpha\in\Delta^-\), and suppose
that \(E^\rho(-a\alpha)\) is nonzero and of type~\(\mathrm{II}\).
The dominant representative of \(\rho-a\alpha\) uniquely
determines \(a\) and \(\alpha\).
Consequently, these \(E^\rho(-a\alpha)\) are linearly independent.
\end{cor}

\begin{pf}
In the two cases of type~\(\mathrm{II}\) in
Lemma~\ref{nonzero E criteria for B complex}, both new coordinates
exceed \(n-\tfrac12\). The assertions follow from
Lemma~\ref{lem:signed-permutation-orbit-recovery}(2).
\end{pf}

Define \(X,Y\subset\ZZ_{\ge1}^2\) by
\begin{align*}
X&=\{(a,b):2\mid b,\ 2\nmid a,\ a\ge2n+1,\
1\le|a-b|\le2n-1\},\\
Y&=\{(a,b):2\mid b,\ 2\nmid a,\ a\le2n-1,\
a+b\ge2n+1,\ 1\le|a-b|\le2n-1\}.
\end{align*}
For \(a,b\in\ZZ_{\ge1}\), put
\[
E_{a,b}^{\mathrm{II},-}
:=\sum_{\alpha\in\Omega_{a,b}^{\mathrm{II}}\cap\Delta^-}
E^\rho(-a\alpha).
\]
By Lemma~\ref{lem:signed-permutation-orbit-recovery}(3),
the span of the \(E^\rho(-a\alpha)\) in
Corollary~\ref{independent type II B} has zero intersection with
\(\operatorname{span}\{A_{a,b}:(a,b)\in\Gamma\}\).
Together with \eqref{B-term-decomposition}, we have
\begin{equation}\label{B-type-split-zero}
M_{p,D}=0
\quad\Longleftrightarrow\quad
E_P^{\mathrm I}=0
\ \text{and}\
E_P^{\mathrm{II}}=0.
\end{equation}

For \((a,b)\in\Gamma\) with \(ab=pD\), define
\(\epsilon_P(a,b)\) by the following table:
\begin{equation}\label{eq:B-epsilon-table}
\begin{array}{c|c|c}
\text{condition on }(a,b)&\epsilon_P(a,b)=1&\epsilon_P(a,b)=-1\\ \hline
2\mid a,\ b\le2n-1
&\text{never}&p\mid2b,\ p\nmid b\\[1mm]
\begin{gathered}2\nmid a,\ a\le2n-1\\b\ne2a\end{gathered}
&p\mid2a,\ p\nmid a&\text{never}\\[1mm]
\text{otherwise}
&p\mid b,\ p\nmid a&p\mid a,\ p\nmid b
\end{array}
\end{equation}
In each row the coefficient is zero unless one of the displayed
divisibility conditions holds. Set \(\epsilon_P(a,b)=0\) also
when \((a,b)\notin\Gamma\) or \(ab\ne pD\).

\begin{prop}\label{B_cr}
We have
\begin{equation}\label{eq:B-exact-coefficient}
M_{p,D}
=\sum_{\substack{(a,b)\in\Gamma\\ab=pD}}
\epsilon_P(a,b)A_{a,b}
+\sum_{(a,b)\in P\cap S_{\mathrm{II},-}}
E_{a,b}^{\mathrm{II},-}
-\sum_{(a,b)\in t(P)'\cap S_{\mathrm{II},-}}
E_{a,b}^{\mathrm{II},-}.
\end{equation}
Moreover, \(M_{p,D}=0\) if and only if
\(
\epsilon_P(a,b)=0\)
for every \((a,b)\in\Gamma\) and \(
P\cap S_{\mathrm{II},-}
=t(P)'\cap S_{\mathrm{II},-}\).
\end{prop}
\begin{pf}
By the definitions of the sets, \(S_{\mathrm I,-}^{(1)}=\Gamma\) and
\[
\begin{aligned}
t(S_{\mathrm I,-}^{(2)})
&=\{(a,b)\in\Gamma:2\mid a,\ b\le2n-1\},\\
S_{\mathrm I,+}^{(3)}\cup S_{\mathrm I,-}^{(4)}
&=\{(a,b)\in\Gamma:2\nmid a,\ a\le2n-1,\ b\ne2a\}.
\end{aligned}
\]
The last two sets are disjoint. In case~(a) of the table in
Lemma~\ref{nonzero term set}, the negative root and its opposite
contribute \(+1\) when \(p\mid b\) and \(-1\) when \(p\mid a\),
respectively. On the first of these two subsets, the roots in
case~(b) and their opposites give additional contributions
\(+1\) when \(p\mid a\) and \(-1\) when \(p\mid2b\).
On the second subset, the roots in cases~(c) and~(d) and their
opposites give additional contributions \(+1\) when \(p\mid2a\)
and \(-1\) when \(p\mid b\).
There are no further contributions. Cancelling the common terms
and using \(p\mid a\Rightarrow p\mid2a\) and
\(p\mid b\Rightarrow p\mid2b\) gives
\eqref{eq:B-epsilon-table}, hence the first sum in
\eqref{eq:B-exact-coefficient}. Transposition and \(t\) preserve
\(ab\), so every resulting pair satisfies \(ab=pD\).
By Corollary~\ref{independent type I B},
\(E_P^{\mathrm I}=0\) if and only if \(\epsilon_P=0\) on \(\Gamma\).

If a positive root \(\alpha\) contributes a term of
type~\(\mathrm{II}\) at \((a,b)\in P\), then \(-\alpha\)
contributes at \((b/2,2a)\), with the opposite sign.
This gives the two sums involving \(E_{a,b}^{\mathrm{II},-}\)
in \eqref{eq:B-exact-coefficient}.
Since \(b=2a-2(\rho|\alpha)\) for long roots,
distinct pairs \((a,b)\) involve disjoint sets of \((a,\alpha)\).
Each \(E_{a,b}^{\mathrm{II},-}\) with
\((a,b)\in S_{\mathrm{II},-}\) contains at least one term.
Corollary~\ref{independent type II B} therefore shows that the
\(E_{a,b}^{\mathrm{II},-}\), with \((a,b)\in S_{\mathrm{II},-}\),
are nonzero and linearly independent. Hence
\[
E_P^{\mathrm{II}}=0
\quad\Longleftrightarrow\quad
P\cap S_{\mathrm{II},-}
=t(P)'\cap S_{\mathrm{II},-}.
\]
Together with \eqref{B-type-split-zero}, this proves the proposition.

\end{pf}

\begin{cor}\label{cor:B-parity}
\begin{enumerate}
\item
If \(p\) is odd, then \(M_{p,D}=0\) if and only if
\begin{align}
U\cap X&=\emptyset,\label{oddCondition1}\\
U\cap S_{\mathrm{II},-}&=\emptyset.\label{oddCondition3}
\end{align}
\item
If \(p\equiv2\pmod4\), then \(M_{p,D}=0\) if and only if
\begin{align}
P\cap X&=t(P)'\cap Y=\emptyset,\label{2Condition1}\\
P\cap S_{\mathrm{II},-}
&=t(P)'\cap S_{\mathrm{II},-}.\label{2Condition3}
\end{align}
\item
If \(4\mid p\), then \(M_{p,D}=0\) if and only if
\begin{align}
P\cap X&=\emptyset,\label{4Condition1}\\
P\cap S_{\mathrm{II},-}
&=t(P)'\cap S_{\mathrm{II},-}.\label{4Condition3}
\end{align}
\end{enumerate}
\end{cor}
\begin{pf}
Suppose first that \(p\) is odd. For the pairs in the first two
rows of \eqref{eq:B-epsilon-table}, we have \(\epsilon_P(a,b)=0\),
since multiplication by \(2\) does not change divisibility by \(p\).
In the last row, \(\epsilon_P(a,b)=0\) also when \(b=2a\).
For every other pair in that row, either \((a,b)\) or \((b,a)\)
belongs to \(X\); conversely, arranging the coordinates of a pair
in \(X\) in increasing order gives a pair in the last row.
Here \(\epsilon_P(a,b)\ne0\) precisely when \(p\) divides exactly
one coordinate. Since \(U=P\mathbin\triangle P'=U'\), this
proves \eqref{oddCondition1}.

Now suppose that \(p\) is even. For a pair in \(\Gamma\) with
\(ab=pD\), write \(o\) for its odd coordinate and \(e\) for
its even coordinate. If \(o\ge2n+1\), then \((o,e)\in X\),
and the last row gives a nonzero coefficient exactly when
\(p\mid e\), equivalently \((o,e)\in P\cap X\).
If \(o\le2n-1\), then \((o,e)\in Y\).
For \(e<o\) the first row applies; for \(o<e\ne2o\) the
second row applies. In both cases the coefficient is nonzero
exactly when \(p\mid2o\), since \(p\nmid o\).
For \(e=2o\), the last row gives the same condition.
As \(e\) is even, \(p\mid2o\) is equivalent to
\((o,e)\in t(P)'\cap Y\). Hence the type~\(\mathrm I\)
coefficients all vanish precisely when
\(P\cap X=t(P)'\cap Y=\emptyset\).
This proves \eqref{2Condition1}. If \(4\mid p\), an odd \(o\)
cannot satisfy \(p\mid2o\), so \(t(P)'\cap Y=\emptyset\)
automatically, giving \eqref{4Condition1}.

For type~\(\mathrm{II}\), Proposition~\ref{B_cr} gives
\eqref{2Condition3} and \eqref{4Condition3}.
If \(p\) is odd, then \(p\mid2a\) if and only if \(p\mid a\).
Since the second coordinate is even on \(S_{\mathrm{II},-}\),
\[
t(P)'\cap S_{\mathrm{II},-}
=P'\cap S_{\mathrm{II},-}.
\]
Thus \(E_P^{\mathrm{II}}=0\) is equivalent to
\(U\cap S_{\mathrm{II},-}=\emptyset\), proving
\eqref{oddCondition3}.
\end{pf}

We now determine   $\lambda_{\nu}=\kappa\Lambda_0-qD_p\delta+\nu$, \(\nu\in\mathcal S^B_{p,n}\), together with the expansion of \(M_{p,D_p}\).

\begin{theorem}\label{Bconclusion}
Let $\mathfrak g$ be of type $B_n$ with $n\ge2$, and assume that 
\[
2(\kappa+\mathbf{h}^\vee)=\frac pq,
\qquad p\in\ZZ_{\ge2},\quad q\in\ZZ_{\ge1},\quad (p,q)=1.
\]
	Let  $\lambda_{\nu}=\kappa\Lambda_0-qD_p\delta+\nu$, \(\nu\in\mathcal S^B_{p,n}\),  be the weight of a minimal singular vector in $V^{\kappa}(\mathfrak{g})$. Then we have 
\begin{enumerate}
\item If $2\mid p$ and $p\ge4n-2$, then
\[
D_p=\frac p2-2n+2,\qquad
\nu=D_p(e_1+e_2).
\]
\item If $2\nmid p$ and $p\ge2n$, then
\[
D_p=p-2n+1,\qquad
\nu=(p-2n+1)e_1.
\]
\item Suppose either that \(p\) is odd, \(5\le p<2n\),
\(\lfloor2n/p\rfloor\ne1\), and \((p,n)\ne(5,9)\), or that
\(p\) is even, \(10\le p<4n-2\), and
\((p,n)\notin\{(12,5),(12,8),(16,7)\}\).
Put
\[
s_1=\left\lfloor\frac{2n}{p}\right\rfloor,\qquad
s_2=s_1+1,\qquad
s_3=\left\lfloor\frac{4n}{p}\right\rfloor,\qquad
s_4=s_3+1 .
\]
Set
\[
D^{(1)}=
\begin{cases}
2s_2(s_2p-2n+2), & 2\nmid p,\ 2\nmid s_1,\\
s_2(s_2p-2n+1), & 2\nmid p,\ 2\mid s_1,\\
s_4\left(\dfrac{s_4p}{2}-2n+2\right), & 2\mid p,
\end{cases}
\]
\[
D^{(2)}=
\begin{cases}
2s_1(2n+2-s_1p),
& 2\nmid p,\ s_1\ge1,\ 2s_1p\ge2n+3,\\
s_3\left(2n+2-\dfrac{s_3p}{2}\right),
& 2\mid p,\ 4\mid s_3p,\ s_3p\ge2n+3,\\
\dfrac{s_3}{2}\left(2n+1-\dfrac{s_3p}{2}\right),
& 2\mid p,\ 4\nmid s_3p,\\
+\infty, & \text{otherwise},
\end{cases}
\]
and
\[
D^{(3)}=
\begin{cases}
\dfrac12\left(2n+1-\dfrac{p}{2}\right),
& p\equiv2\pmod4,\ s_3\ge1,\ 2\mid s_3,\\
+\infty, & \text{otherwise}.
\end{cases}
\]
Set
\[
\nu^{(1)}=
\begin{cases}
(s_2p-2n+2)(e_1+e_2), & 2\nmid p,\ 2\nmid s_1,\\
(s_2p-2n+1)e_1, & 2\nmid p,\ 2\mid s_1,\\
\left(\dfrac{s_4p}{2}-2n+2\right)(e_1+e_2), & 2\mid p.
\end{cases}
\]
If $D^{(2)}<+\infty$, set
\[
\nu^{(2)}=
\begin{cases}
2(e_1+\cdots+e_{2n+2-s_1p}), & 2\nmid p,\\
2(e_1+\cdots+e_{2n+2-\frac{s_3p}{2}}), & 2\mid p,\ 4\mid s_3p,\\
e_1+\cdots+e_{2n+1-\frac{s_3p}{2}}, & 2\mid p,\ 4\nmid s_3p,
\end{cases}
\]
and if  $D^{(3)}<+\infty$, 
\[
\nu^{(3)}=e_1+\cdots+e_{p/2}.
\]
Then
\[
D_p=\min\{D^{(1)},D^{(2)},D^{(3)}\},
\]
and \(\mathcal S^B_{p,n}\) consists of the $\nu=\nu^{(i)}$ for which \(D^{(i)}=D_p\). 
\end{enumerate}
 Finally, we have
\begin{equation}\label{eq:B-first-coefficient}
	M_{p,D_p}=\sum_{\nu\in\mathcal S^B_{p,n}}E^\rho(\nu).
\end{equation}
The minimal conformal weight of \(N_\kappa(\mathfrak g)\) is
\(qD_p\).
\end{theorem}
\begin{theorem}\label{specialcaseB}
Let \(n\ge2\) and
\[
2(\kappa+\mathbf{h}^\vee)=\frac pq,\qquad
p\ge2,\quad q\ge1,\quad (p,q)=1,
\]
and put \(s_1=\lfloor2n/p\rfloor\). The cases not covered by
Theorem~\ref{Bconclusion} are as follows. 
\begin{enumerate}
\item We have
\[
\begin{array}{c|c|c}
\text{conditions} & D_p & \nu\\ \hline
p=2 & n & e_1\\[1mm]
p=4 & 2n+1 & 3e_1+e_2\\[1mm]
(p,n)=(5,9) & 32 & 4(e_1+e_2)\ \text{and}\ 3e_1+2e_2+2e_3+2e_4+e_5\\[1mm]
(p,n)=(12,5) & 7 & 3e_1+2e_2+2e_3+2e_4+2e_5\\[1mm]
(p,n)=(12,8) & 12 & 4(e_1+e_2)\\[1mm]
(p,n)=(16,7) & 8 & 4(e_1+e_2).
\end{array}
\]
\item If $p=3$, then
\[
D_p=\begin{cases}
\dfrac{8n}{3},&n\equiv0\pmod3,\\
4\lfloor n/3\rfloor+2,&n\equiv1\pmod3,\\
4n,&n\equiv2\pmod3,
\end{cases}
\]
and
\[
\nu=\begin{cases}
2(e_1+e_2),&n\equiv0\pmod3,\\
2e_1,&n\equiv1\pmod3,\\
3e_1+2e_2+e_3,&n\equiv2\pmod3,\ n\ne2,\\
3e_1+2e_2,&n=2.
\end{cases}
\]
\item If $p=6$ and $n\ge3$, then
\[
D_p=\begin{cases}
\dfrac{2n-1}{3},&n\equiv2\pmod3,\\
n-1,&n\not\equiv2\pmod3,
\end{cases}
\]
and
\[
\nu=\begin{cases}
e_1+e_2,&n\equiv2\pmod3,\\
e_1+e_2+e_3,&n\not\equiv2\pmod3.
\end{cases}
\]
\item If $p=8$ and $n\ge3$, then
\[
D_p=\begin{cases}n,&2\mid n,\\2n-1,&2\nmid n,\end{cases}
\]
and
\[
\nu=\begin{cases}
2(e_1+e_2),&2\mid n,\\
3e_1+2e_2+2e_3+e_4,&2\nmid n,\ n>3,\\
3e_1+2e_2+2e_3,&n=3.
\end{cases}
\]
\item Suppose $2\nmid p$, $p\ge5$, and $s_1=1$. If $2n+2\ne2p$, then
\[
D_p=\min\{4(2p-2n+2),\,4n+4-2p\},
\]
and if \(D_p=4(2p-2n+2)\), 
\[
\nu=(2p-2n+2)(e_1+e_2).\]
If \(D_p=4n+4-2p\), then
\[\nu=2(e_1+\cdots+e_{2n+2-p}).
\]
If $2n+2=2p$, then
\[
D_p=\min\{16,\,2p+2\},
\]
where  \(\nu=4(e_1+e_2)\) when
\(D_p=16\), and \(\nu=3e_1+2e_2+\cdots+2e_n\) when \(D_p=2p+2\).
\end{enumerate}
\end{theorem}
We prove Theorem \ref{Bconclusion} and Theorem  \ref{specialcaseB} by the following  Lemma \ref{lem:B-endpoint} and Proposition \ref{prop:B-first-layer}.  The assertion about minimal conformal weights follows from Proposition~\ref{prop:first-jantzen-layer}.

 Recall from (\ref{ce1}) that for each \((a,b)\in\ZZ_{\ge1}^2\) with \(p\mid ab\), we denote
\(D=ab/p\) and \(P=P_{p,D}\). 
\begin{lem}\label{lem:B-endpoint} 
	Suppose that \((p,n)\) satisfies the hypotheses of one of the
	three cases in Theorem~\ref{Bconclusion}. If \(M_{p,D}\ne0\), then
	\[
	D\ge
	\begin{cases}
		p/2-2n+2,&2\mid p,\quad p\ge4n-2,\\
		p-2n+1,&2\nmid p,\quad p\ge2n,\\
		\min\{D^{(1)},D^{(2)},D^{(3)}\},&\text{under the hypotheses of case~(3).}
	\end{cases}
	\]
	Here \(D^{(1)},D^{(2)},D^{(3)}\) are defined in
	Theorem~\ref{Bconclusion}.
\end{lem}

\begin{pf}
	By \eqref{eq:B-exact-coefficient}, the coefficient of
	\(E_{a,b}^{\mathrm{II},-}\), for \((a,b)\in S_{\mathrm{II},-}\)
	with \(ab=pD\), is \(1\) when \(p\mid b\) and \(p\nmid2a\),
	and \(-1\) when \(p\mid2a\) and \(p\nmid b\); it is zero otherwise.
	The set \(S_{\mathrm{II},-}\) consists of the pairs with \(2\mid b\)
	satisfying either of the equivalent bounds
	\begin{align}
		\max\{b-4n+4,\,4n+4-b,\,2\}
		&\le2a\le b-2,\label{eq:B-II-a-range}\\
		\max\{2a+2,\,4n+4-2a\}
		&\le b\le2a+4n-4.\label{eq:B-II-b-range}
	\end{align}
	For fixed \(b\), \(D=ab/p\) increases with \(a\); for fixed
	\(a\), it increases with \(b\). The least values allowed by these
	bounds are as follows, where \(k\ge1\) and \(2\mid kp\):
	\begin{align}
		b=kp,\quad 2n+4\le kp\le4n:
		&\quad a=2n+2-kp/2,\quad D=k(2n+2-kp/2),\label{eq:B-branch-b-lower}\\
		b=kp,\quad kp\ge4n+2:
		&\quad a=kp/2-2n+2,\quad D=k(kp/2-2n+2),\label{eq:B-branch-b-upper}\\
		2a=kp,\quad 4\le kp\le2n:
		&\quad b=4n+4-kp,\quad D=k(2n+2-kp/2),\label{eq:B-branch-a-lower}\\
		2a=kp,\quad kp\ge2n+2:
		&\quad b=kp+2,\quad D=k(kp/2+1).\label{eq:B-branch-a-upper}
	\end{align}
	The value of \(D\) strictly decreases as \(k\) increases in
	\eqref{eq:B-branch-b-lower}, and strictly increases in the other
	three ranges. On \(2a+b=4n+4\), with \(p\mid b\) or \(p\mid2a\),
	the coefficient is zero if and only if \(p\mid4n+4\).
	On \(b-2a=4n-4\), with \(p\mid b\), it is zero if and only if
	\(p\mid4n-4\). On \(b=2a+2\), with \(p\mid2a\), it is zero
	if and only if \(p=2\). In particular, all coefficients of
	type~\(\mathrm{II}\) vanish for \(p=2\).
	
	For type~\(\mathrm I\), we use the sets \(X,Y\) and
	Corollary~\ref{cor:B-parity}. If \(2\nmid p\), a pair in
	\(P'\cap X\) has \(a=kp\), with \(2\nmid k\) and \(kp\ge2n+1\),
	and gives \(D\ge k(kp-2n+1)\). A pair in \(P\cap X\) has
	\(b=kp\), with \(2\mid kp\), and gives
	\(D\ge k\max\{2n+1,kp-2n+1\}\). Both bounds strictly increase
	with \(k\). For \(2\nmid p\), such a pair contributes only when
	\(p\) divides exactly one coordinate.
	
	For \(2\mid p\), the set \(t(P)'\cap Y\) is empty unless
	\(p\equiv2\pmod4\) and \(p\le4n-2\). Under these conditions,
	its pairs have \(a=kp/2\), with \(2\nmid k\), \(kp/2\le2n-1\),
	and \(b\ge2n+1-kp/2\). Write \(2n=s_3p/2+r\), with
	\(0\le r<p/2\). The least value of \(k(2n+1-kp/2)/2\) is
	\(s_3(r+1)/2\), attained only at \(k=s_3\), when \(2\nmid s_3\);
	it is \((2n+1-p/2)/2\), attained only at \(k=1\), when \(2\mid s_3\).
	
	We also need the nonzero coefficients when \(p\mid4n+4\).
	Suppose \(p\ge5\). For fixed \(2a=kp\le2n\), the least
	\(b\) giving a nonzero coefficient is \(4n+6-kp\); for fixed
	\(b=kp\) with \(2n+4\le kp\le4n\), the least such \(a\) is
	\(2n+3-kp/2\). Both give \(D=k(2n+3-kp/2)\).
	Together with \eqref{eq:B-branch-a-upper}, these bounds show that,
	among pairs with \(p\mid2a\), or with \(p\mid b\) and \(b\le4n\),
	the least nonzero value is \(4n+6-2p\) at
	\((p,4n+6-2p)\) when \(2\nmid p\). When \(2\mid p\) and
	\(p\le2n+2\), it is \(2n+3-p/2\) at
	\((p/2,4n+6-p)\). These pairs belong to \(S_{\mathrm{II},-}\)
	and have coefficient \(-1\).
	
	\smallskip
	
	\noindent\textbf{Case~1.} \(n\ge2,\ 2\mid p,\ p\ge4n-2\).
	For \(p=4n-2\) and \(p=4n\), the bounds above give \(D\ge1\)
	and \(D\ge2\), respectively. For \(p>4n\), the least value in
	\eqref{eq:B-branch-b-upper} is \(p/2-2n+2\). The least values
	in \(P\cap X\) and \eqref{eq:B-branch-a-upper} are
	\(p-2n+1\) and \(p/2+1\), both larger; the other ranges are empty.
	Thus \(D\ge p/2-2n+2\).
	
	\smallskip
	
	\noindent\textbf{Case~2.} \(n\ge2,\ 2\nmid p,\ p>2n\).
	The least value from \(P'\cap X\) is \(p-2n+1\).
	The least values in \eqref{eq:B-branch-b-upper},
	\eqref{eq:B-branch-a-upper}, and \(P\cap X\) are
	\(2(p-2n+2)\), \(2(p+1)\), and \(2(2p-2n+1)\), respectively.
	All are larger, and the other ranges are empty. Hence
	\(M_{p,D}=0\) for \(0<D<p-2n+1\).
	
	\smallskip
	
	\noindent\textbf{Case~3.} \(n\ge2,\ 2\nmid p,\ 5\le p\le n,\ (p,n)\ne(5,9)\).
	Write \(2n=s_1p+r\), with \(0\le r<p\); then \(s_1\ge2\)
	and \(r\equiv s_1\pmod2\).
	For \(b\le4n\), the least value with \(p\mid b\) occurs at
	\((r+2,2s_1p)\) and is \(2s_1(r+2)=D^{(2)}\).
	For \(b>4n\), it occurs at \((p-r+2,2s_2p)\) and is
	\(2s_2(p-r+2)\).
	For fixed \(a=kp\le n\), the least value is
	\(2k(2n+2-kp)\), which is at least \(D^{(2)}\) for \(k\ge2\).
	For \(k=1\), it can be smaller only when \(r=p-1\); then
	\(2\mid s_1\) and it is larger than \(D^{(1)}=2s_2\).
	For \(n<kp\le2n\), the least value \(2k(kp+1)\) is larger
	than \(D^{(2)}\); for \(kp>2n\), it is larger than
	\(2s_2(p-r+2)\).
	
	The coefficient at \((r+2,2s_1p)\) is zero only when \(r=p-2\).
	Here \(2\nmid s_1\), and the preceding calculation for
	\(p\mid4n+4\) gives the least nonzero value
	\(2s_1p+2\), attained at \((p,2s_1p+2)\).
	It is larger than \(8s_2=D^{(1)}\), except for
	\((p,n)=(5,9)\), when both are \(32\).
	The coefficient at \((p-r+2,2s_2p)\) is zero only when \(r=2\).
	In that case \(2\mid s_1\), the pair \((4,2s_1p)\) has
	nonzero coefficient, and its value \(D^{(2)}=8s_1\) is smaller
	than \(2s_2p\). Thus these cancellations cause no further case.
	
	For \(2\mid s_1\), the least value from \(P'\cap X\) is
	\(s_2(p-r+1)=D^{(1)}\), at \((s_2p,p-r+1)\). Its coefficient
	is nonzero since \(2\mid r\). For \(2\nmid s_1\), the least
	value from \(P'\cap X\) is \((s_1+2)(2p-r+1)\).
	It exceeds \(2s_2(p-r+2)\) when \(r\ge3\), and exceeds
	\(D^{(2)}=6s_1\) when \(r=1\).
	The least value from \(P\cap X\) is \(2(2n+1)\).
	It is larger than \(D^{(1)}\) when \(2\mid s_1\), larger
	than \(D^{(2)}\) when \(2\nmid s_1\) and \(r\le p-4\),
	and larger than \(D^{(1)}=8s_2\) when \(r=p-2\).
	These bounds give \(D\ge\min\{D^{(1)},D^{(2)}\}\), with
	\(D^{(3)}=+\infty\).
	
	\smallskip
	
	\noindent\textbf{Case~4.} \(n\ge2,\ 2\mid p,\ 10\le p<4n-2\),
	with \((p,n)\notin\{(12,5),(12,8),(16,7)\}\).
	Write \(2n=s_3p/2+r\), with \(0\le r<p/2\).
	For \(p\mid b\) and \(b>4n\), the least value is
	\(s_4(p/2-r+2)=D^{(1)}\), at \((p/2-r+2,s_4p)\).
	For \(p\mid b\) and \(b\le4n\), it is \(s_3(r+2)\), at
	\((r+2,s_3p)\), provided \(s_3p\ge2n+4\).
	For fixed \(2a=kp\le2n\), the least value
	\(k(2n+2-kp/2)\) is at least \(s_3(r+2)\) for \(k\ge2\).
	For \(k=1\), a smaller value can occur only when \(r=p/2-1\).
	Then \(p\equiv2\pmod4\), \(2\mid s_3\), and this value is
	larger than \(D^{(3)}=s_3p/4\).
	For \(2n<kp\le4n\), the least value \(k(kp/2+1)\) is larger
	than \(s_3(r+2)\); for \(kp>4n\), it is larger than \(D^{(1)}\).
	
	The coefficient at \((r+2,s_3p)\) is zero when \(r=p/2-2\).
	If \(4\nmid s_3p\), then \(p\equiv2\pmod4\), \(2\nmid s_3\),
	and the type~\(\mathrm I\) value \(D^{(2)}=s_3(r+1)/2\) is
	smaller than \(s_3p/2\). If \(4\mid s_3p\), then \(4\mid p\).
	For \(s_3\ge2\), the inequality \(s_3p/2\le4s_4\) leaves only
	\((p,n)=(12,8)\); all other such parameters have
	\(D^{(1)}=4s_4<s_3p/2\).
	For \(s_3=1\), the pair \((p/2,p)\) is not in
	\(S_{\mathrm{II},-}\). The least nonzero values with
	\(a=p/2\) fixed, with \(b=2p\) fixed, and in \(P\cap X\)
	are \(p/2+1\), \(8\), and \(p-1\), respectively.
	The parameters \(p=12,16\) give \((12,5),(16,7)\); for
	\(p\ge20\), the least value is \(8\).
	When \(s_3\ge2\), the calculation for \(p\mid4n+4\) gives
	\(D\ge s_3p/2+1\) for the nonzero terms in the other ranges.
	
	The coefficient at \((p/2-r+2,s_4p)\) is zero only when \(r=2\).
	Then \((4,s_3p)\in S_{\mathrm{II},-}\) has nonzero coefficient
	and gives \(D=4s_3<D^{(1)}\); the bounds with \(2a=kp\)
	also give \(D\ge4s_3\).
	Finally, if \(p\equiv2\pmod4\), the least value from
	\(t(P)'\cap Y\) is \(s_3(r+1)/2=D^{(2)}\) when \(2\nmid s_3\),
	and \((2n+1-p/2)/2=D^{(3)}\) when \(2\mid s_3\).
	For \(4\mid s_3p\), the type~\(\mathrm{II}\) value
	\(s_3(r+2)\), when its pair belongs to \(S_{\mathrm{II},-}\),
	is \(D^{(2)}\). The least value \(2n+1\) from \(P\cap X\)
	is larger than \(\min\{D^{(1)},D^{(2)},D^{(3)}\}\).
	This proves the required bound in Case~4 and completes the proof.
\end{pf}

For fixed \(p,n\), let \(D_0\) be the value of the formula for
\(D_p\) in Theorem~\ref{Bconclusion} or Theorem~\ref{specialcaseB}.
Let \(\mathcal S_0\) be the corresponding set of weights, with
\(D_0\) in place of \(D_p\).

\begin{prop}\label{prop:B-first-layer}
	We have
	\[
	M_{p,D}=0\qquad,
	\]
for $0<D<D_0$	and
	\begin{equation}\label{eq:B-remaining-first-layer}
		M_{p,D_0}
		=\sum_{\nu\in\mathcal S_0}E^\rho(\nu).
	\end{equation}
	The weights in \(\mathcal S_0\) are pairwise distinct. In particular,
	\(D_p=D_0\) and \(\mathcal S^B_{p,n}=\mathcal S_0\).
\end{prop}

\begin{pf}
	We use \eqref{eq:B-exact-coefficient} and the bounds in
	Lemma~\ref{lem:B-endpoint}.
	
	\smallskip
	
	\noindent\textbf{Case~1.} \(n\ge2\), with either \(2\mid p,\ p\ge4n-2\), or \(2\nmid p,\ p>2n\).
	Lemma~\ref{lem:B-endpoint} gives \(M_{p,D}=0\) for \(0<D<D_0\).
	For \(2\nmid p\), the only pair contributing to \(M_{p,D_0}\)
	is \((p-2n+1,p)\in\Gamma\); its coefficient and Weyl sign are
	both \(1\), giving \(E^\rho((p-2n+1)e_1)\).
	For \(p=4n-2\), only \((2,p/2)\in\Gamma\) contributes, with
	\(\epsilon_P(2,p/2)=-1\) and \(A_{2,p/2}=-E^\rho(e_1+e_2)\).
	For \(2\mid p\) and \(p\ge4n\), the only pair is
	\((a,b)=(p/2-2n+2,p)\in S_{\mathrm{II},-}\).
	The coefficient of \(E_{a,b}^{\mathrm{II},-}\) is \(1\), and
	\(\Omega_{a,b}^{\mathrm{II}}\cap\Delta^-\) contains only
	\(-(e_1+e_2)\), giving \(E^\rho((p/2-2n+2)(e_1+e_2))\).
	
	\smallskip
	
	\noindent\textbf{Case~2.} \(n\ge2\), with either
	\(2\nmid p,\ 5\le p\le n,\ (p,n)\ne(5,9)\), or
	\(2\mid p,\ 10\le p<4n-2,\ (p,n)\notin\{(12,5),(12,8),(16,7)\}\).
	Lemma~\ref{lem:B-endpoint} gives \(M_{p,D}=0\) for \(0<D<D_0\).
	Write \(2n=s_1p+r\), with \(0\le r<p\), when \(2\nmid p\),
	and \(2n=s_3p/2+r\), with \(0\le r<p/2\), when \(2\mid p\).
	Equality in the preceding bounds gives the following pairs of
	type~\(\mathrm{II}\):
	\[
	\begin{array}{c|c|c}
		\text{condition}&(a,b)&D\\ \hline
		2\nmid p,\ 2\nmid s_1&(p-r+2,2s_2p)&D^{(1)}\\
		2\nmid p&(r+2,2s_1p)&D^{(2)}\\
		2\mid p&(p/2-r+2,s_4p)&D^{(1)}\\
		2\mid p,\ 4\mid s_3p&(r+2,s_3p)&D^{(2)}.
	\end{array}
	\]
	The last row requires \(s_3p\ge2n+4\).
	A pair in the table contributes to \(M_{p,D_0}\) precisely when
	\(ab=pD_0\) and \(p\nmid2a\); its coefficient is then \(1\).
	For \(2\nmid p\) and \(2\mid s_1\), the smaller value
	\(D^{(1)}=s_2(p-r+1)\) comes from type~\(\mathrm I\).
	For \(2\mid p\) and \(4\nmid s_3p\), the smaller value
	\(D^{(2)}=s_3(r+1)/2\) also comes from type~\(\mathrm I\).
	
	The type~\(\mathrm I\) coefficients and weights are as follows.
	For \(2\nmid p\) and \(2\mid s_1\),
	\(\epsilon_P(p-r+1,s_2p)=1\) and
	\(A_{p-r+1,s_2p}=E^\rho((p-r+1)e_1)\).
	For \(2\mid p\) and \(4\nmid s_3p\),
	\(\epsilon_P(r+1,s_3p/2)=-1\) and
	\(A_{r+1,s_3p/2}=-E^\rho(e_1+\cdots+e_{r+1})\).
	For \(p\equiv2\pmod4\) and \(2\mid s_3\),
	\(\epsilon_P(p/2,2n+1-p/2)=1\) and
	\(A_{p/2,2n+1-p/2}=E^\rho(e_1+\cdots+e_{p/2})\).
	Each pair is used only when \(ab=pD_0\).
	The bounds on its coordinates ensure that all weight indices are
	at most \(n\).
	
	For \((a,b)\in S_{\mathrm{II},-}\) with \(b=2a+4n-4\),
	the only root in \(\Omega_{a,b}^{\mathrm{II}}\cap\Delta^-\) is
	\(\alpha=-(e_1+e_2)\), and
	\(E^\rho(-a\alpha)=E^\rho(a(e_1+e_2))\).
	For \(b=4n+4-2a\), it is \(\alpha=-(e_{a-1}+e_a)\), and
	\(E^\rho(-a\alpha)=E^\rho(2(e_1+\cdots+e_a))\).
	Thus every pair with \(ab=pD_0\) and a nonzero coefficient in
	the preceding lists contributes \(E^\rho(\nu)\) with coefficient \(1\).
	The equality conditions in the bounds exclude other pairs.
	The listed pairs are distinct, so
	Lemma~\ref{lem:signed-permutation-orbit-recovery} gives distinct
	highest weights and linearly independent \(E^\rho(\nu)\).
	
	\smallskip
	
	\noindent\textbf{Case~3.} \(n\ge2,\ p\in\{2,3,4\}\), or \(n\ge3,\ p\in\{6,8\}\).
	The contributions to \(M_{p,D_0}\) are listed in Table~\ref{tab:B-small-p}. In each row of type~\(\mathrm I\),
	\(\epsilon_P(a,b)=\sigma\) and
	\(A_{a,b}=\sigma E^\rho(\nu)\).
	In each row of type~\(\mathrm{II}\), the coefficient multiplying
	\(E_{a,b}^{\mathrm{II},-}\) in \eqref{eq:B-exact-coefficient}
	is \(\sigma\), and
	\(E_{a,b}^{\mathrm{II},-}=\sigma E^\rho(\nu)\).
	\begin{center}\captionof{table}{Contributions to $M_{p,D_0}$ in type $B_n$
for $p\in\{2,3,4,6,8\}$.}
\label{tab:B-small-p}
		\small
		\begin{adjustbox}{max width=\textwidth}
			\begin{tabular}{c|c|c|c|c|l}
				conditions & type & \((a,b)\) & \(D\) & \(\sigma\) & \(\nu\)\\ \hline
				\(p=2\)&\(\mathrm I\)&\((1,2n)\)&\(n\)&\(1\)&\(e_1\)\\
				\(p=4\)&\(\mathrm I\)&\((4,2n+1)\)&\(2n+1\)&\(-1\)&\(3e_1+e_2\)\\
				\(p=3,\ n\equiv0\pmod3\)&\(\mathrm{II}\)&\((2,4n)\)&\(8n/3\)&\(1\)&\(2(e_1+e_2)\)\\
				\(p=3,\ n\equiv1\pmod3\)&\(\mathrm I\)&\((2,2n+1)\)&\((4n+2)/3\)&\(1\)&\(2e_1\)\\
				\(p=3,\ n\equiv2\pmod3,\ n>2\)&\(\mathrm{II}\)&\((3,4n)\)&\(4n\)&\(-1\)&\(3e_1+2e_2+e_3\)\\
				\(p=3,\ n=2\)&\(\mathrm{II}\)&\((3,8)\)&\(8\)&\(-1\)&\(3e_1+2e_2\)\\
				\(p=6,\ n\not\equiv2\pmod3\)&\(\mathrm I\)&\((3,2n-2)\)&\(n-1\)&\(1\)&\(e_1+e_2+e_3\)\\
				\(p=6,\ n\equiv2\pmod3\)&\(\mathrm I\)&\((2,2n-1)\)&\((2n-1)/3\)&\(-1\)&\(e_1+e_2\)\\
				\(p=8,\ 2\mid n\)&\(\mathrm{II}\)&\((2,4n)\)&\(n\)&\(1\)&\(2(e_1+e_2)\)\\
				\(p=8,\ 2\nmid n,\ n>3\)&\(\mathrm{II}\)&\((4,4n-2)\)&\(2n-1\)&\(-1\)&\(3e_1+2e_2+2e_3+e_4\)\\
				\(p=8,\ n=3\)&\(\mathrm{II}\)&\((4,10)\)&\(5\)&\(-1\)&\(3e_1+2e_2+2e_3\)
			\end{tabular}
		\end{adjustbox}
	\end{center}
	For each type~\(\mathrm{II}\) pair in the table, the set
	\(\Omega_{a,b}^{\mathrm{II}}\cap\Delta^-\) contains exactly one root:
	\(-(e_1+e_2)\) at \((2,4n)\),
	\(-(e_1+e_3)\) at \((3,4n)\) for \(n>2\),
	\(-(e_1-e_2)\) at \((3,8)\) for \(n=2\),
	\(-(e_2+e_4)\) at \((4,4n-2)\) for \(2\nmid n\) and \(n>3\), and
	\(-(e_2-e_3)\) at \((4,10)\) for \(n=3\).
	For each of these roots, \(E^\rho(-a\alpha)=\sigma E^\rho(\nu)\),
	with \(\sigma\) and \(\nu\) as in the table.
	
	For \(p=2\), all type~\(\mathrm{II}\) coefficients vanish.
	The minimum on \(t(P)'\cap Y\) is \(n\), whereas the lower
	bound on \(P\cap X\) is \(2n+1\).
	For \(p=4\), the set \(t(P)'\cap Y\) is empty and the minimum
	on \(P\cap X\) is \(2n+1\). A nonzero type~\(\mathrm{II}\)
	coefficient requires \(2\mid a,\ b\equiv2\pmod4\), or
	\(2\nmid a,\ 4\mid b\). For \(a=2\), \eqref{eq:B-II-a-range} gives
	\(b=4n\), whose type~\(\mathrm{II}\) coefficient is zero.
	For \(3\le a\le n+1\), the bound \(b\ge4n+6-2a\) gives
	\(ab\ge12n\); for \(a\ge n+1\), the same conclusion follows
	from \(b\ge2a+2\). Equality is attained at \((3,4n)\), giving
	\(D=3n>2n+1\).
	
	For \(p=3\), the least values of \(D\) for which
	\(E_{P_{p,D}}^{\mathrm I}\ne0\) and
	\(E_{P_{p,D}}^{\mathrm{II}}\ne0\), respectively, are
	\[
	\begin{array}{c|c|c}
		n&\min\{D:E_{P_{p,D}}^{\mathrm I}\ne0\}&
		\min\{D:E_{P_{p,D}}^{\mathrm{II}}\ne0\}\\ \hline
		3m&8m+4&8m\\
		3m+1&4m+2&12m+2\\
		3m+2&12m+10&12m+8.
	\end{array}
	\]
	The first two rows have \(m\ge1\), and the last has \(m\ge0\).
	For type~\(\mathrm I\), the three values follow from the bounds
	on \(X\) at the beginning of the proof of Lemma~\ref{lem:B-endpoint}.
	For type~\(\mathrm{II}\), take \(p=3\) in
	\eqref{eq:B-branch-b-lower}--\eqref{eq:B-branch-a-upper}.
	If \(n=3m\), the pair \((2,4n)\) gives \(D=8m\).
	If \(n=3m+1\), the pair \((3,4n-2)\) gives \(D=12m+2\).
	If \(n=3m+2>2\), the coefficient at \((3,4n-2)\) is zero,
	since \(3\mid2n+2\). With \(a=3\) fixed, the least larger
	even \(b\) with nonzero coefficient is \(4n\), giving
	\(D=12m+8\). For \(n=2\), \((3,6)\notin S_{\mathrm{II},-}\)
	and \((3,8)\) gives \(D=8\). In the last case the least
	\(D\) with \(E_{P_{p,D}}^{\mathrm I}\ne0\) is \(4n+2\).
	
	For \(p=8\), \(2\nmid n\), and \(n>3\), the coefficient at
	\((4,4n-4)\) is zero.
	With \(a=4\) fixed, the least larger even \(b\) giving a
	nonzero coefficient is \(4n-2\).
	For \(n=3\), \((4,8)\notin S_{\mathrm{II},-}\).
	The least \(D\) with a nonzero type~\(\mathrm I\) coefficient
	is \(2n+1\). The coefficients and highest weights in the table
	follow by substitution in \eqref{eq:B-exact-coefficient} and
	Lemma~\ref{nonzero E criteria for B complex}.
	For \(p=6\), the type~\(\mathrm{II}\) divisibility conditions
	are the same as for \(p=3\), while \(D=ab/6\).
	The corresponding minima are therefore half the entries in the
	last column above. The minimum on \(t(P)'\cap Y\) is \(n-1\) when
	\(n\not\equiv2\pmod3\), and \((2n-1)/3\) otherwise; each is
	strictly below the type~\(\mathrm{II}\) minimum and \(2n+1\),
	the minimum on \(P\cap X\).
	For \(p=8\), the set \(t(P)'\cap Y\) is empty. If \(2\mid n\), \((a,b)=(2,4n)\) gives \(D=n\); the other three cases in
	\eqref{eq:B-branch-b-lower}--\eqref{eq:B-branch-a-upper} give larger
	values of \(D\). If \(2\nmid n\), then \(8\mid4n+4\), and
	the preceding calculation gives \(D=2n-1\) at \((4,4n-2)\).
	Both values are less than \(2n+1\). Thus, for
	\(p\in\{2,3,4,6,8\}\), we have \(M_{p,D}=0\) for
	\(0<D<D_0\), and the table gives every term of \(M_{p,D_0}\).
	
	\smallskip
	
	\noindent\textbf{Case~4.} \((p,n)\in\{(5,9),(12,5),(12,8),(16,7)\}\).
	The type~\(\mathrm{II}\) contributions to \(M_{p,D_0}\) are listed in Table~\ref{tab:B-exceptional},
	with the same convention for \(\sigma\):
	\begin{center}
    \captionof{table}{Type~$\mathrm{II}$ contributions to $M_{p,D_0}$
in type $B_n$ for
$(p,n)\in\{(5,9),(12,5),(12,8),(16,7)\}$.}
\label{tab:B-exceptional}
		\small
		\begin{adjustbox}{max width=\textwidth}
			\begin{tabular}{c|c|c|c|l}
				\((p,n)\)&\((a,b)\)&\(D\)&\(\sigma\)&\(\nu\)\\ \hline
				\((5,9)\)&\((4,40)\)&\(32\)&\(1\)&\(4(e_1+e_2)\)\\
				\((5,9)\)&\((5,32)\)&\(32\)&\(-1\)&\(3e_1+2e_2+2e_3+2e_4+e_5\)\\
				\((12,5)\)&\((6,14)\)&\(7\)&\(-1\)&\(3e_1+2e_2+2e_3+2e_4+2e_5\)\\
				\((12,8)\)&\((4,36)\)&\(12\)&\(1\)&\(4(e_1+e_2)\)\\
				\((16,7)\)&\((4,32)\)&\(8\)&\(1\)&\(4(e_1+e_2)\)
			\end{tabular}
		\end{adjustbox}
	\end{center}
	For \((p,n)=(12,5)\), the only root in
	\(\Omega_{6,14}^{\mathrm{II}}\cap\Delta^-\) is \(-(e_4-e_5)\).
	For \((p,n)=(12,8)\) or \((16,7)\), the corresponding set
	\(\Omega_{a,b}^{\mathrm{II}}\cap\Delta^-\) contains only
	\(-(e_1+e_2)\). For \((p,n)=(5,9)\), the type~\(\mathrm{II}\)
	coefficient at \((5,30)\) is zero since \(5\mid30\) and \(5\mid10\).
	The pairs \((a,b)=(4,40)\) and
	\((5,32)\), both with \(D=32\), have respectively
	\(-(e_1+e_2)\) and \(-(e_3+e_5)\) as the only roots in
	\(\Omega_{a,b}^{\mathrm{II}}\cap\Delta^-\). The least \(D\)
	with a nonzero type~\(\mathrm I\) coefficient is \(38\).
	For \((p,n)\in\{(5,9),(12,5),(12,8),(16,7)\}\), the
	comparisons in the proof of Lemma~\ref{lem:B-endpoint} give
	\(M_{p,D}=0\) for \(0<D<D_0\), and the table gives all terms of
	\(M_{p,D_0}\). In particular, both terms for \((p,n)=(5,9)\)
	occur with coefficient \(1\).
	
	\smallskip
	
	\noindent\textbf{Case~5.} \(n\ge2,\ 2\nmid p,\ p\ge5,\ n<p\le2n\).
	Write \(2n=p+r\), where \(1\le r\le p-2\) and \(2\nmid r\).
	For \(r\le p-4\), the pairs \((p-r+2,4p)\) and \((r+2,2p)\)
	give \(D=4(p-r+2)\) and \(D=2(r+2)\), respectively.
	Both have type~\(\mathrm{II}\) coefficient \(1\), and their
	contributions are \(E^\rho((p-r+2)(e_1+e_2))\) and
	\(E^\rho(2(e_1+\cdots+e_{r+2}))\).
	For fixed \(a=p\), the least value is \(2p+2>2(r+2)\).
	The least values from \(P\cap X\) and \(P'\cap X\) are
	\(2(2n+1)>2(r+2)\) and \(3(2p-r+1)>4(p-r+2)\), respectively.
	Consequently, \(D_0=\min\{4(p-r+2),2(r+2)\}\), and only the
	pair giving this value contributes. The two values are distinct:
	the first is divisible by \(4\), whereas the second is congruent
	to \(2\) modulo \(4\).
	
	For \(r=p-2\), the pair \((p,2p)\) is not in
	\(S_{\mathrm{II},-}\). The least nonzero values come from
	\((4,4p)\) and \((p,2p+2)\), giving \(D=16\) and \(D=2p+2\).
	Their coefficients are \(1\) and \(-1\), respectively, and
	\(E_{4,4p}^{\mathrm{II},-}=E^\rho(4(e_1+e_2))\), while
	\(E_{p,2p+2}^{\mathrm{II},-}
	=-E^\rho(3e_1+2e_2+\cdots+2e_n)\).
	The least values from type~\(\mathrm I\) are \(4p-2\) and
	\(3p+9\), both larger than \(\min\{16,2p+2\}\).
	Thus \(D_0=\min\{16,2p+2\}\). Both pairs contribute only for
	\((p,n)=(7,6)\), when their highest weights are distinct.
	Together with Cases~1--4 and
	Lemma~\ref{lem:signed-permutation-orbit-recovery}, these calculations
	give \eqref{eq:B-remaining-first-layer} and the asserted distinctness
	of the highest weights.
\end{pf}

\section{\texorpdfstring{Type $C_n$}{Type Cn}}\label{sec:type-C}
Assume \(n\ge3\). We use the standard coordinate realization of the root system of type $C_n$ \cite{Humphreys,tauvel2005lie} as
\[
\Delta_{\mathrm{short}}=\{\pm(e_i\pm e_j):1\leq i<j\leq n\},\qquad
\Delta_{\mathrm{long}}=\{\pm2e_i:1\leq i\leq n\},
\]
where \((e_i|e_j)=\tfrac12\delta_{ij}\). The positive roots are
\[
\Delta^+=\{e_i\pm e_j:1\leq i<j\leq n\}\cup\{2e_i:1\leq i\leq n\}.
\]
We take \(\alpha_i=e_i-e_{i+1}\) for \(1\le i<n\) and
\(\alpha_n=2e_n\). The Weyl vector is \(\rho=(n,n-1,\dots,1)\). For fixed
\(p\) and \(D\), put \(P=P_{p,D}\)(see (\ref{ce1})). We use
\(\omega_i=e_1+\cdots+e_i\) for the fundamental weights of \(C_n\).

For \(a\in\ZZ_{\ge1}\) and \(\alpha\in\Delta\) with
\(E^\rho(-a\alpha)\ne0\), let \(\mu\) be the dominant representative
of \(\rho-a\alpha\).
We say that \(E^\rho(-a\alpha)\) is of type~\(\mathrm I\)
if \(\mu\) has exactly one coordinate greater than \(n\),
and of type~\(\mathrm{II}\) if it has exactly two such coordinates.
\begin{lem}\label{nonzero E criteria for C complex}
Let \(a\in\ZZ_{\ge1}\) and \(\alpha\in\Delta^-\).
Then \(E^\rho(-a\alpha)\ne0\) precisely in the following cases.
Here \(1\le i\le n\) for \(\alpha=-2e_i\), and
\(1\le i<j\le n\) for roots involving two indices.
\begin{enumerate}[(i)]
	\item Type $\mathrm{I}$:
	\begin{enumerate}[(a)]
\item \(\alpha=-2e_i\) and \(2a\ge i\).
\item \(\alpha=-(e_i+e_j)\) and \(a=j-i\ge i\).
\item \(\alpha=-(e_i-e_j)\) and \(a=2n-i-j+2\ge i\).
\item \(\alpha=-(e_i-e_j)\) and \(a=2n-2j+2\ge i\).
\end{enumerate}
The dominant representative in case~(a) is
\[
(n-i+1+2a,n,\dots,\widehat{n-i+1},\dots,1).
\]
In cases~(b) and~(c) it is
\[
(n-i+1+a,n,\dots,\widehat{n-j+1},\dots,1),
\]
and in case~(d) it is
\[
(n-i+1+a,n,\dots,\widehat{n-i+1},\dots,1).
\]

	\item Type $\mathrm{II}$:
	\begin{enumerate}[(a)]
		\item $\alpha=-(e_i+e_j)$ and $a\geq j$, with dominant representative
		\[
		(n-i+1+a,n-j+1+a,n,\dots,\widehat{n-i+1},\dots,\widehat{n-j+1},\dots,1).
		\]
		\item $\alpha=-(e_i-e_j)$ and $a\geq 2n-j+2$, with dominant representative
		\[
		(n-i+1+a,a-(n-j+1),n,\dots,\widehat{n-i+1},\dots,\widehat{n-j+1},\dots,1).
		\]
	\end{enumerate}
\end{enumerate}

\end{lem}
\begin{pf}
Apply Lemma~\ref{lem:signed-permutation-alternants} to the three
families of negative roots. For \(\alpha=-2e_i\), the changed coordinate
is equal to an unchanged coordinate when \(2a<i\), and exceeds
\(n\) when \(2a\ge i\). This gives case~(i)(a).

For \(\alpha=-(e_i+e_j)\), regularity requires \(a\ge i\).
If \(a<j\), the second changed coordinate must equal the omitted
\(i\)-th coordinate, so \(a=j-i\). Otherwise both changed
coordinates exceed \(n\). These are cases~(i)(b) and~(ii)(a).

For \(\alpha=-(e_i-e_j)\), again \(a\ge i\) is necessary.
The second changed coordinate is equal to an unchanged coordinate
when \(a<n-j+1\), and is zero when \(a=n-j+1\).
For larger \(a\), its absolute value \(a-(n-j+1)\) must equal one
of the omitted coordinates \(n-i+1,n-j+1\), or exceed \(n\).
These alternatives give \(a=2n-i-j+2\), \(a=2n-2j+2\), and
\(a\ge2n-j+2\), respectively. Sorting gives exactly the
representatives in the remaining cases.
\end{pf}
For \(\tau\in\{\mathrm{I},\mathrm{II}\}\), define
\[
\Omega_{a,b}^{\tau}
:=\{\alpha\in\Omega_{a,b}:
E^\rho(-a\alpha)\ne0
\text{ and }E^\rho(-a\alpha)\text{ is of type }\tau\},
\]
and put
\[
E_{a,b}^{\tau}
:=\sum_{\alpha\in\Omega_{a,b}^{\tau}}E^\rho(-a\alpha),
\qquad
E_P^{\tau}
:=\sum_{(a,b)\in P}E_{a,b}^{\tau}.
\]
Every nonzero term in \(E_{a,b}\) belongs to exactly one of these
two types. Hence
\begin{equation}\label{C-term-decomposition}
M_{p,D}=E_P^{\mathrm{I}}
+E_P^{\mathrm{II}}.
\end{equation}
The nonzero terms described in
Lemma~\ref{nonzero E criteria for C complex} are indexed by the following
subsets of \(\ZZ_{\geq1}^2\):
\begin{center}
\small
\(\displaystyle
\begin{aligned}
&S_{\mathrm{I},-}^{(1)}:=
\left\{(a,b):
\begin{gathered}
2\mid b,\ 2\leq b-2a\leq 2n,\\
2n+2\leq 2a+b
\end{gathered}
\right\},\\
&S_{\mathrm{I},-}^{(2)}:=
\left\{(a,b):
\begin{gathered}
2\mid b,\ 2+2a\leq b\leq 2n,\\
2n+2\leq 2a+b
\end{gathered}
\right\},\\
&S_{\mathrm{I},-}^{(3)}:=
\left\{(a,b):
\begin{gathered}
2\mid b,\ a+1\leq b\leq 2a-2,\\
b\leq 2n,\ 2n+2\leq 2a+b
\end{gathered}
\right\},\\
&S_{\mathrm{I},-}^{(4)}:=
\left\{(a,b):
\begin{gathered}
2\mid a,\ 2b-2n\leq a\leq b-1,\\
a\leq 2n,\ 2n+2\leq 2b+a
\end{gathered}
\right\},\\
&S_{\mathrm{II},-}:=
\left\{(a,b):
\begin{gathered}
2\leq a,\ 1\leq b-a\leq 2n-1,\\
b+a\geq 2n+3
\end{gathered}
\right\}
\setminus\{(n+1,n+2),(n+1,n+3)\}.
\end{aligned}
\)
\end{center}
Define
\[
S_{\mathrm{I},+}^{(1)}=t(S_{\mathrm{I},-}^{(1)})^{\prime}, \ \
S_{\mathrm{I},+}^{(j)}=(S_{\mathrm{I},-}^{(j)})^{\prime}
\quad (j=2,3,4), \ \ 
S_{\mathrm{II},+}=(S_{\mathrm{II},-})^{\prime}.
\]
Set
\[
S_{\mathrm{I},\pm}=\bigcup_j S_{\mathrm{I},\pm}^{(j)},\ \ 
S_{\mathrm{I}}=S_{\mathrm{I},-}\cup S_{\mathrm{I},+},\ \
S_{\mathrm{II}}=S_{\mathrm{II},-}\cup S_{\mathrm{II},+}.
\]
The sets \(S_{\mathrm I}\) and \(S_{\mathrm{II}}\) need not be disjoint. For example, in type \(C_4\),
\(E_{3,10}\) contains the type~\(\mathrm{I}\) term associated with
\(\alpha=-2e_3\) and the type~\(\mathrm{II}\) term associated with
\(\alpha=-(e_1+e_2)\).

Put \(\Gamma=S_{\mathrm I,-}^{(1)}\).
For \((a,b)\in\Gamma\), put
\[
A_{a,b}
:=(-1)^{n-\frac b2+a}
E\!\left(
\frac b2+a,n,\dots,
\widehat{\frac b2-a},\dots,1
\right).
\]
Put \(j=n+1-\tfrac b2+a\in\{1,\ldots,n\}\). Then
\[
\nu=\left(\frac b2+a-n\right)e_1+\sum_{i=2}^{j}e_i,
\qquad
A_{a,b}=(-1)^{j-1}E^\rho(\nu),
\]
where the sum is empty when \(j=1\).

\begin{lem}\label{EPI=0 for C}
For \(a,b\in\ZZ_{\ge1}\) and
\(\tau\in\{\mathrm I,\mathrm{II}\}\), the sum \(E_{a,b}\)
contains a nonzero term of type \(\tau\) if and only if
\((a,b)\in S_\tau\).
\end{lem}
\begin{pf}
For type~\(\mathrm I\), substituting the four cases for negative
roots in Lemma~\ref{nonzero E criteria for C complex}
into the equations defining \(\Omega_{a,b}\) gives
\[
\begin{array}{c|c|c|c}
\text{case}&i&j&E^\rho(-a\alpha)\\ \hline
\mathrm{(a)}&n-\frac b2+a+1&-&A_{a,b}\\
\mathrm{(b)}&n-\frac b2+1&n-\frac b2+a+1&-A_{a,b}\\
\mathrm{(c)}&n-\frac b2+1&n+\frac b2-a+1&A_{b/2,2a}\\
\mathrm{(d)}&n-b+\frac a2+1&n-\frac a2+1&-A_{a/2,2b}.
\end{array}
\]
Under these substitutions, the conditions on \(i,j,a\) in
cases~(a)--(d) are equivalent to those defining
\(S_{\mathrm{I},-}^{(1)},\ldots,S_{\mathrm{I},-}^{(4)}\),
respectively.
The signs follow from the coordinate permutations and sign
changes taking \(\rho-a\alpha\) to its dominant representative.

The term for \(-\alpha\) has the opposite sign and occurs at
\((b/2,2a)\) in case~(a), or at \((b,a)\) in cases~(b)--(d).
This gives \(S_{\mathrm{I},+}^{(1)},\ldots,S_{\mathrm{I},+}^{(4)}\).

For type~\(\mathrm{II}\) and \(\alpha\in\Delta^-\),
Lemma~\ref{nonzero E criteria for C complex} and the equations
defining \(\Omega_{a,b}\) give
\[
b=
\begin{cases}
a+2n-i-j+2,
&\alpha=-(e_i+e_j),\quad a\ge j,\\
a+j-i,
&\alpha=-(e_i-e_j),\quad a\ge2n-j+2.
\end{cases}
\]
In the first case, \(3\le b-a\le2n-1\) and
\(a+b\ge2n+3\). In the second, \(1\le b-a\le n-1\) and
\(a+b\ge2n+5\). These inequalities exclude
\((n+1,n+2)\) and \((n+1,n+3)\), so
\((a,b)\in S_{\mathrm{II},-}\).

Conversely, let \((a,b)\in S_{\mathrm{II},-}\) and put
\(c=b-a\). Then \(1\le c\le2n-1\) and
\(2a+c\ge2n+3\). If \(c\ge3\), take
\[
j=\left\lceil n+1-\frac{c-1}{2}\right\rceil,
\qquad i=2n+2-c-j,
\qquad \alpha=-(e_i+e_j).
\]
The bounds on \(c\) give \(1\le i<j\le n\), and
\(2a+c\ge2n+3\) gives \(a\ge j\).
Since \(2(\rho|\alpha)=-c=a-b\), we have
\(\alpha\in\Omega_{a,b}\), and
Lemma~\ref{nonzero E criteria for C complex} gives a nonzero
term of type~\(\mathrm{II}\).

If \(c=1\), take \(\alpha=-(e_{n-1}-e_n)\); if \(c=2\),
take \(\alpha=-(e_{n-2}-e_n)\).
The excluded pairs \((n+1,n+2)\) and \((n+1,n+3)\), together
with \(2a+c\ge2n+3\), give \(a\ge n+2\).
In both cases, \(2(\rho|\alpha)=-c=a-b\), and the same
lemma gives a nonzero term of type~\(\mathrm{II}\).

All roots contributing terms of type~\(\mathrm{II}\) are short.
Replacing \(\alpha\) by \(-\alpha\) sends \((a,b)\) to
\((b,a)\) and reverses the sign of the corresponding term.
Hence \(S_{\mathrm{II},+}=(S_{\mathrm{II},-})'\).
\end{pf}
\begin{cor}\label{independent type I C}
The elements \(A_{a,b}\), with \((a,b)\in\Gamma\), are nonzero
and linearly independent.
\end{cor}

\begin{pf}
Put \(u=a+b/2\) and \(v=b/2-a\). Since \((a,b)\in\Gamma\),
we have \(u\ge n+1\) and \(v\in\{1,\ldots,n\}\).
The coordinates of the dominant representative are \(u\) and the
coordinates of \(\rho\) other than \(v\). Since \(a=(u-v)/2\) and \(b=u+v\), the assertion
follows from Lemma~\ref{lem:signed-permutation-orbit-recovery}(1).
\end{pf}
\begin{cor}\label{EPII=0 for C}
Let \(a\in\ZZ_{\ge1}\) and \(\alpha\in\Delta^-\), and suppose
that \(E^\rho(-a\alpha)\) is nonzero and of type~\(\mathrm{II}\).
The dominant representative of \(\rho-a\alpha\) uniquely
determines \(a,\alpha\).
Consequently, these \(E^\rho(-a\alpha)\) are linearly independent.
\end{cor}

\begin{pf}
In the two cases of type~\(\mathrm{II}\) in
Lemma~\ref{nonzero E criteria for C complex}, both new coordinates
exceed \(n\). The assertions follow from
Lemma~\ref{lem:signed-permutation-orbit-recovery}(2).
\end{pf}
For \(a,b\in\ZZ_{\ge1}\), put
\[
E_{a,b}^{\mathrm{II},-}
:=\sum_{\alpha\in\Omega_{a,b}^{\mathrm{II}}\cap\Delta^-}
E^\rho(-a\alpha).
\]
By Lemma~\ref{lem:signed-permutation-orbit-recovery}(3),
the span of the \(E^\rho(-a\alpha)\) in
Corollary~\ref{EPII=0 for C} has zero intersection with
\(\operatorname{span}\{A_{a,b}:(a,b)\in\Gamma\}\).

Together with \eqref{C-term-decomposition}, this gives
\begin{equation}\label{C-type-split-zero}
M_{p,D}=0
\quad\Longleftrightarrow\quad
E_P^{\mathrm{I}}=0
\ \text{and}\
E_P^{\mathrm{II}}=0.
\end{equation}

Put \(\Gamma(pD)=\{(a,b)\in\Gamma:ab=pD\}\), and partition
\(\Gamma\) into the four regions
\[
\begin{aligned}
 \mathcal R_1&=\{(a,b)\in\Gamma:b\le2n,\ b\ne4a\},\\
 \mathcal R_2&=\{(a,b)\in\Gamma:b\le2n,\ b=4a\},\\
 \mathcal R_3&=\{(a,b)\in\Gamma:b\ge2n+2,\ a\le n,\ b\ne4a\},\\
 \mathcal R_4&=\{(a,b)\in\Gamma:b\ge2n+2,\
                 (a\ge n+1\text{ or }b=4a)\}.
\end{aligned}
\]
For \((a,b)\in\Gamma(pD)\), define \(\epsilon_P(a,b)\) by
the following table, with value zero in every case not listed:
\begin{equation}\label{eq:C-pointwise-delta}
\begin{array}{c|c|c}
\text{region}&\epsilon_P(a,b)=1&\epsilon_P(a,b)=-1\\ \hline
\mathcal R_1&p\mid a,\ 2p\nmid b&p\nmid a,\ 2p\mid b\\
\mathcal R_2&\text{never}&p\mid2a,\ p\nmid a\\
\mathcal R_3&p\mid b,\ 2p\nmid b&\text{never}\\
\mathcal R_4&p\mid b,\ p\nmid2a&p\nmid b,\ p\mid2a
\end{array}
\end{equation}
For \((a,b)\notin\Gamma(pD)\), set \(\epsilon_P(a,b)=0\).

\begin{prop}\label{C_cr}
We have
\begin{equation}\label{eq:C-exact-coefficient}
M_{p,D}
={}\sum_{(a,b)\in\Gamma(pD)}\epsilon_P(a,b)A_{a,b}+\sum_{(a,b)\in P\cap S_{\mathrm{II},-}}
 E_{a,b}^{\mathrm{II},-}
-\sum_{(a,b)\in P'\cap S_{\mathrm{II},-}}
 E_{a,b}^{\mathrm{II},-}.
\end{equation}
Moreover, \(M_{p,D}=0\) if and only if
\[
\epsilon_P|_{\Gamma(pD)}=0
\quad\text{and}\quad
U\cap S_{\mathrm{II},-}=\emptyset.
\]
\end{prop}
\begin{pf}
The changes of indices in Lemma~\ref{EPI=0 for C} preserve \(ab\).
Fix \(\gamma=(a,b)\in\Gamma(pD)\) and put
\(\widehat\gamma=(b/2,2a)\).
The long roots contribute \(+1\) when \(p\mid b\) and \(-1\)
when \(p\mid2a\) to the coefficient of \(A_{a,b}\).
For the short roots, the conditions in
Lemma~\ref{EPI=0 for C} give
\[
\begin{aligned}
\gamma\in S_{\mathrm{I},-}^{(2)}
&\Longleftrightarrow b\le2n,\\
\widehat\gamma\in S_{\mathrm{I},-}^{(3)}
&\Longleftrightarrow a\le n,\ b\le4a-2,\\
\widehat\gamma\in S_{\mathrm{I},+}^{(4)}
&\Longleftrightarrow a\le n,\ b\ge4a+2.
\end{aligned}
\]
For \(b\le2n\), the contributions from the negative root in
case~(b) and its opposite are \(-1\) when \(p\mid b\) and \(+1\)
when \(p\mid a\), respectively. Since \(b\) is even, one of the
last two conditions holds precisely when \(a\le n\) and \(b\ne4a\).
The roots in cases~(c) and~(d) and their opposites then contribute
\(+1\) when \(p\mid2a\) and \(-1\) when \(2p\mid b\).

On \(\mathcal R_1\), both short-root contributions occur; after
cancellation, the remaining contributions are \(+1\) when
\(p\mid a\) and \(-1\) when \(2p\mid b\).
On \(\mathcal R_2\), only case~(b) occurs, leaving \(+1\) when
\(p\mid a\) and \(-1\) when \(p\mid2a\).
On \(\mathcal R_3\), case~(b) is absent, leaving \(+1\) when
\(p\mid b\) and \(-1\) when \(2p\mid b\).
On \(\mathcal R_4\), neither short-root contribution occurs, so
only the long-root contributions remain. Combining opposite signs
in each region gives exactly \eqref{eq:C-pointwise-delta}.

All roots contributing terms of type~\(\mathrm{II}\) are short.
Replacing a positive root by its negative sends \((a,b)\) to
\((b,a)\) and reverses the sign, giving the two sums involving
\(E_{a,b}^{\mathrm{II},-}\) in \eqref{eq:C-exact-coefficient}.
Since \(b=a-2(\rho|\alpha)\) for short roots,
distinct pairs \((a,b)\) involve disjoint sets of \((a,\alpha)\).
Each \(E_{a,b}^{\mathrm{II},-}\) with
\((a,b)\in S_{\mathrm{II},-}\) contains at least one term.
Corollary~\ref{EPII=0 for C} shows that the
\(E_{a,b}^{\mathrm{II},-}\), with
\((a,b)\in S_{\mathrm{II},-}\), are nonzero and linearly independent.
Hence \(E_P^{\mathrm{II}}=0\) if and
only if \(P\cap S_{\mathrm{II},-}=P'\cap S_{\mathrm{II},-}\).
By Corollary~\ref{independent type I C},
\(E_P^{\mathrm I}=0\) if and only if \(\epsilon_P=0\)
on \(\Gamma(pD)\).
Together with \eqref{C-type-split-zero}, this proves the proposition.
\end{pf}
\begin{cor}
\label{lem:C-pointwise-coefficients}
If \(p\) is odd, then \(E_P^{\mathrm I}=0\) if and only if
\[
 (P\mathbin\triangle P')\cap
 (\mathcal R_1\cup\mathcal R_4)=\emptyset.
\]
If \(p\) is even, then \(E_P^{\mathrm I}=0\) if and only if every
\((a,b)\in\Gamma(pD)\) satisfies the condition for its region:
\[
\begin{array}{c|c}
\mathcal R_1&p\mid a\ \Longleftrightarrow\ 2p\mid b\\
\mathcal R_2&p\mid2a\ \Longrightarrow\ p\mid a\\
\mathcal R_3&p\mid b\ \Longrightarrow\ 2p\mid b\\
\mathcal R_4&p\mid b\ \Longleftrightarrow\ p\mid2a
\end{array}
\]
\end{cor}
\begin{pf}
If \(p\) is odd, then \(p\mid2a\) is equivalent to \(p\mid a\).
Since \(b\) is even, \(2p\mid b\) is equivalent to \(p\mid b\).
Thus \(\epsilon_P(a,b)=0\) on \(\mathcal R_2\cup\mathcal R_3\).
On each of the other two regions, \(\epsilon_P(a,b)\ne0\) precisely
when \(p\) divides exactly one of \(a,b\). This proves the odd case.
For even \(p\), the four conditions are equivalent to
\(\epsilon_P(a,b)=0\) in the respective regions, by
\eqref{eq:C-pointwise-delta}. The conclusion follows from
Corollary~\ref{independent type I C}.
\end{pf}

We now determine $\lambda_{\nu}=\kappa\Lambda_0-qD_p\delta+\nu$, \(\nu\in\mathcal S^C_{p,n}\), together with the expansion of \(M_{p,D_p}\).
\begin{theorem}\label{C_conclusion}
Let $\mathfrak g$ be of type $C_n$ with \(n\ge3\), and assume
\[
2(\kappa+\mathbf{h}^\vee)=\frac pq,
\qquad p\in\ZZ_{\ge5},\quad q\in\ZZ_{\ge1},\quad (p,q)=1.
\]
Exclude \((p,n)\in\{(5,4),(5,9),(6,11)\}\), which are treated in
Theorem~\ref{specialcaseC}. Put
\[
s_1=\left\lfloor\frac{2n}{p}\right\rfloor,\qquad s_2=s_1+1.
\]
As in Theorem~\ref{Bconclusion}, the value \(+\infty\) is omitted
from a minimum. Let  $\lambda_{\nu}=\kappa\Lambda_0-qD_p\delta+\nu$,$ \nu\in \mathcal S^C_{p,n}$,   be the weight of a minimal singular vector in $V^{\kappa}(\mathfrak{g})$. 
Then:
\begin{enumerate}
\item If $2\mid p$ and $p\ge2n+2$, then
\[
D_p=\frac p2-n,\qquad \nu=(p-2n)e_1.
\]
\item If $2\nmid p$ and $p\ge2n$, then
\[
D_p=p-2n+1,\qquad \nu=(p-2n+1)(e_1+e_2).
\]
\item Suppose that \(s_1=1\). If \(p\ne n+2\), then
\[
D_p=\min\{2n+3-p,\,4p-4n+2\}.
\]
\(\nu=2(e_1+\cdots+e_{2n+3-p})\) when \(D_p=2n+3-p\), in which
case \(2n+3-p\le n\), and \(\nu=(2p-2n+1)(e_1+e_2)\) when
\(D_p=4p-4n+2\).

If  $p=n+2$, then
$
D_p=\min\{p+1,10\},
$
where 
\(\nu=4e_1+2(e_2+\cdots+e_n)\) when \(D_p=p+1\), and
\(\nu=5(e_1+e_2)\) when \(D_p=10\).
\item Suppose that \(s_1\ge2\). Define
\[
D^{(1)}=\begin{cases}
s_2(s_2p-2n+1),&2\nmid p,\ 2\mid s_1,\\
s_2\left(\dfrac{s_2p}{2}-n\right),&2\mid p,\ 2\mid s_1,\\
(s_1-1)\left(n+1-\dfrac{(s_1-1)p}{2}\right),&2\nmid s_1,
\end{cases}
\]
\[
D^{(2)}=\begin{cases}
s_1\left(n-\dfrac{s_1p}{2}+1\right),&2\mid s_1,\\
s_2(s_2p-2n+1),&2\nmid s_1,
\end{cases}
\qquad
D^{(3)}=\begin{cases}+\infty,&2\mid s_1,\\
s_1(2n+3-s_1p),&2\nmid s_1
\end{cases}
\]
Then
\[
D_p=\min\{D^{(1)},D^{(2)},D^{(3)}\}.
\]
If \(D^{(1)}=D_p\), set
\[
\nu^{(1)}=\begin{cases}
(s_2p-2n+1)(e_1+e_2),&2\nmid p,\ 2\mid s_1,\\
(s_2p-2n)e_1,&2\mid p,\ 2\mid s_1,\\
e_1+\cdots+e_{2n+2-(s_1-1)p},&2\nmid s_1
\end{cases}
\]
If \(2\nmid s_1\) and \(D^{(1)}=D_p\), then
\(2n+2-(s_1-1)p\le n\). If \(D^{(2)}=D_p\), set
\[
\nu^{(2)}=\begin{cases}
e_1+\cdots+e_{2n+2-s_1p},&2\mid s_1,\\
(s_2p-2n+1)(e_1+e_2),&2\nmid s_1
\end{cases}
\]
When $2\mid s_1$, one has $2n+2-s_1p\le n$. Finally, if
$2\nmid s_1$ and $D^{(3)}=D_p$, then
\[
\nu^{(3)}=2(e_1+\cdots+e_{2n+3-s_1p}),
\]
and $2n+3-s_1p\le n$. \(\mathcal S^C_{p,n}\) consists of the
\(\nu^{(i)}\) for which \(D^{(i)}=D_p\).
\end{enumerate}
Finally, we have
\begin{equation}\label{eq:C-first-coefficient}
M_{p,D_p}
=\sum_{\nu\in\mathcal S^C_{p,n}}E^\rho(\nu).
\end{equation}
The minimal conformal weight of \(N_\kappa(\mathfrak g)\) is
\(qD_p\).
\end{theorem}
\begin{theorem}\label{specialcaseC}
Let $\mathfrak g$ be of type $C_n$ with \(n\ge3\), and assume
\[
2(\kappa+\mathbf{h}^\vee)=\frac pq,
\qquad p\in\ZZ_{\ge2},\quad q\in\ZZ_{\ge1},\quad (p,q)=1.
\]
The cases not covered by Theorem~\ref{C_conclusion} are as follows.
In each case, the minimal conformal weight of
\(N_\kappa(\mathfrak g)\) is \(qD_p\).
\begin{enumerate}
\item If $p=2$, then
\[
D_p=\begin{cases}n,&2\mid n,\\2n,&2\nmid n,\end{cases}
\qquad
\nu=\begin{cases}e_1+e_2,&2\mid n,\\2e_1+e_2+e_3,&2\nmid n.\end{cases}
\]
\item If $p=3$, then
\[
D_p=\begin{cases}
\dfrac{2n}{3},&n\equiv0\pmod3,\\
4\lfloor n/3\rfloor,&n\equiv1\pmod3,\\
2n-2,&n\equiv2\pmod3,
\end{cases}
\]
and
\[
\nu=\begin{cases}
e_1+e_2,&n\equiv0\pmod3,\\
e_1+e_2+e_3+e_4,&n\equiv1\pmod3,\\
2e_1+e_2+e_3+e_4+e_5,&n\equiv2\pmod3.
\end{cases}
\]
\item If $p=4$, then
\[
D_p=\begin{cases}
\dfrac n2,&n\equiv0\pmod4,\\
\dfrac{n+1}{2},&n\equiv1\pmod4,\\
\dfrac{3n-6}{2},&n\equiv2\pmod4,\\
10,&n=7,\\
\dfrac{3n+3}{2},&n\equiv3\pmod4,\ n\ne7.
\end{cases}
\]
The finite highest weights are
\[
\begin{array}{c|l}
\text{condition}&\nu\\ \hline
n\equiv0\pmod4&e_1+e_2\\
n\equiv1\pmod4&2e_1\\
n\equiv2\pmod4&e_1+\cdots+e_6\\
n=7&2e_1+e_2+\cdots+e_7\\
n=11&2e_1+e_2+\cdots+e_7,\quad 3e_1+3e_2\\
n\equiv3\pmod4,\ n\ne7,11&3e_1+3e_2.
\end{array}
\]
\item We have
\[
\begin{array}{c|c|c}
(p,n)&D_p&\nu\\ \hline
(5,4)&6&3(e_1+e_2)\\
(5,9)&12&3e_1+3e_2,\quad 2e_1+e_2+\cdots+e_9\\
(6,11)&12&3(e_1+e_2).
\end{array}
\]
\end{enumerate}
Finally, 
\begin{equation}\label{eq:C-remaining-first-coefficient}
M_{p,D_p}
=\sum_{\nu\in\mathcal S^C_{p,n}}E^\rho(\nu).
\end{equation}
\end{theorem}
We prove Theorem \ref{C_conclusion} and  Theorem \ref{specialcaseC} by the following Lemma \ref{lem:C-endpoint} and 
Proposition~\ref{prop:C-first-layer}. 
The assertion about minimal conformal weights follows from
Proposition~\ref{prop:first-jantzen-layer}.
\begin{lem}\label{lem:C-endpoint}
	Suppose that \(n\ge3\) and either \(p>2n\) or \(5\le p\le n\).
	If \(M_{p,D}\ne0\), then
	\[
	D\ge
	\begin{cases}
		p/2-n,&2\mid p,\quad p\ge2n+2,\\
		p-2n+1,&2\nmid p,\quad p>2n,\\
		\min\{D^{(1)},D^{(2)},D^{(3)}\},&5\le p\le n.
	\end{cases}
	\]
\end{lem}
\begin{pf}
	By \eqref{eq:C-exact-coefficient}, the coefficient of
	\(E_{a,b}^{\mathrm{II},-}\), for \((a,b)\in S_{\mathrm{II},-}\)
	with \(ab=pD\), is \(1\) when \(p\mid b\) and \(p\nmid a\),
	and \(-1\) when \(p\mid a\) and \(p\nmid b\); it is zero otherwise.
	The two possibilities for a nonzero coefficient are
	\begin{align}
		b=kp,\quad p\nmid a:
		&\quad \max\{2,kp-2n+1,2n+3-kp\}\le a\le kp-1,
		\quad D=ka,\label{eq:C-II-b-range}\\
		a=kp,\quad p\nmid b:
		&\quad \max\{kp+1,2n+3-kp\}\le b\le kp+2n-1,
		\quad D=kb.\label{eq:C-II-a-range}
	\end{align}
	Here \(k\ge1\), and the pairs \((n+1,n+2),(n+1,n+3)\) are
	excluded. For fixed \(b\), \(D=ab/p\) increases with \(a\);
	for fixed \(a\), it increases with \(b\). In
	\eqref{eq:C-II-b-range}, the least value allowed by the inequalities
	is \(k(2n+3-kp)\) for \(n+2\le kp\le2n+1\), and
	\(k(kp-2n+1)\) for \(kp\ge2n+1\). These values strictly decrease
	and strictly increase with \(k\), respectively.
	
	For type~\(\mathrm I\), a nonzero coefficient requires
	\(p\mid b\) or \(p\mid2a\). If \(b=kp\), then \(2\mid kp\),
	\(D=ka\), and
	\begin{equation}\label{eq:C-I-b-range}
		\max\left\{1,\frac{kp}{2}-n,n+1-\frac{kp}{2}\right\}
		\le a\le\frac{kp}{2}-1.
	\end{equation}
	For fixed \(a\) with \(p\mid2a\), the corresponding bounds are
	\begin{equation}\label{eq:C-I-a-range}
		b\in2\ZZ,\qquad
		\max\{2a+2,2n+2-2a\}\le b\le2a+2n.
	\end{equation}
	For \(n+2\le b=kp\le2n\), the least value is
	\(k(n+1-kp/2)\); for \(b=kp\ge2n+2\), it is
	\(k(kp/2-n)\). For \(2a=kp\), the least value is
	\(k(n+1-kp/2)\) when \(kp\le n\), and \(k(kp/2+1)\)
	when \(kp\ge n+1\). In each case the integer pair must satisfy
	\eqref{eq:C-I-b-range} or \eqref{eq:C-I-a-range}, and its coefficient
	is given by \eqref{eq:C-pointwise-delta}.
	
	\smallskip
	
	\noindent\textbf{Case~1.} \(n\ge3,\ 2\mid p,\ p\ge2n+2\).
	For type~\(\mathrm I\) with \(p\mid b\), the least value is
	\(p/2-n\), at \((p/2-n,p)\). If \(p\nmid b\), then
	\(p\mid2a\), so \(a\ge p/2\) and \(b\ge2a+2\) give
	\(D\ge p/2+1\). For type~\(\mathrm{II}\), the bounds are
	\(D\ge p-2n+1\) when \(p\mid b\), and \(D\ge p+1\)
	when \(p\mid a\). All three are larger than \(p/2-n\).
	
	\smallskip
	
	\noindent\textbf{Case~2.} \(n\ge3,\ 2\nmid p,\ p>2n\).
	For type~\(\mathrm{II}\) with \(p\mid b\), the least value is
	\(p-2n+1\), at \((p-2n+1,p)\). The bound with \(p\mid a\)
	is \(D\ge p+1\). For type~\(\mathrm I\), the conditions
	\(p\mid b\) and \(2\mid b\) give \(b\ge2p\), hence
	\(D\ge2(p-n)\); if \(p\nmid b\), then \(p\mid a\) and
	\(D\ge2p+2\). These values are all larger than \(p-2n+1\).
	
	\smallskip
	
	\noindent\textbf{Case~3.} \(n\ge3,\ 2\nmid p,\ 5\le p\le n\).
	Write \(2n=s_1p+r\), with \(0\le r<p\). On
	\(\mathcal R_2\cup\mathcal R_3\), the coefficient is zero.
	On \(\mathcal R_1\), a nonzero coefficient requires either
	\(b=kp\) with \(2\mid k\), or \(a=mp\).
	For \(2\mid k\) and \(2\le k\le s_1\), the bound
	\(D\ge k(n+1-kp/2)\) is at least \(D^{(2)}\) when
	\(2\mid s_1\), with equality only at \(k=s_1\).
	When \(2\nmid s_1\), it is at least \(D^{(1)}\), with equality
	only at \(k=s_1-1\), or at \(k=2\) if \(r=p-2\).
	For \(a=mp\), \eqref{eq:C-I-a-range} gives the same bound with
	\(k=2m\); equality requires \(b=2n+2-2mp\).
	
	If \(p\mid n+1\), write \(n+1=\ell p\). Then
	\(s_1=2\ell-1\), \(r=p-2\), and \(D^{(2)}=6\ell\).
	For \(a=mp\) in \(\mathcal R_1\), a nonzero coefficient gives
	\(D\ge2m((\ell-m)p+1)\) when \(2m\le\ell\), and
	\(D\ge2m(mp+1)\) when \(2m\ge\ell+1\).
	In the first range, \((mp,2(\ell-m)p)\) has zero coefficient
	or does not belong to \(\Gamma\); the least possible second
	coordinate with a nonzero coefficient is \(2(\ell-m)p+2\).
	These bounds are greater than \(6\ell\), except at
	\((p,n)=(5,9)\), where the least value is \(12\).
	For \(\mathcal R_4\), a nonzero coefficient requires
	\(b\ne4a\) and \(a\ge n+1\). The two divisibility conditions give
	\(D\ge s_2(n+1)\) and \(D\ge s_2(n+2)\), respectively, both
	at least \(s_2(p-r+1)\).
	
	For type~\(\mathrm{II}\) with \(p\mid b\), the least values
	allowed by \eqref{eq:C-II-b-range} are \(s_1(r+3)\) for
	\(b\le2n\), and \(s_2(p-r+1)\) for \(b\ge2n+1\).
	They occur only at \(b=s_1p\) and \(b=s_2p\), respectively,
	when the corresponding pairs belong to \(S_{\mathrm{II},-}\).
	For \(a=kp\) and \(p\nmid n+1\), \eqref{eq:C-II-a-range}
	gives \(D>D^{(2)}\) when \(2\mid s_1\), and
	\(D>D^{(1)}\) when \(2\nmid s_1\).
	For \(n+1=\ell p\), it gives \(D>6\ell\). In particular,
	when \(k=\ell\), the two excluded pairs force
	\(b\ge\ell p+3\), hence \(D\ge\ell(\ell p+3)>6\ell\).
	
	If \(2\mid s_1\), then \(D^{(1)}=s_2(p-r+1)\),
	\(D^{(2)}=s_1(r+2)/2<s_1(r+3)\), and \(D^{(3)}=+\infty\).
	If \(2\nmid s_1\), then
	\(D^{(1)}=(s_1-1)(p+r+2)/2\), \(D^{(2)}=s_2(p-r+1)\),
	and \(D^{(3)}=s_1(r+3)\). Thus every nonzero term satisfies the
	required lower bound.
	
	\smallskip
	
	\noindent\textbf{Case~4.} \(n\ge3,\ 2\mid p,\ 6\le p\le n\).
	Write \(2n=s_1p+r\), with \(0\le r<p\); then \(2\mid r\).
	The estimates in Case~3 for \(\mathcal R_1\) and for
	\((P'\setminus P)\cap S_{\mathrm{II},-}\) use only \(p\ge5\)
	and remain valid. For \(\mathcal R_4\) with \(b\ne4a\),
	we again have \(a\ge n+1\), giving
	\(D\ge s_2(n+1)\) when \(p\mid b\), and
	\(D\ge s_2(n+2)\) when \(p\mid2a\).
	
	On \(\mathcal R_3\), a nonzero coefficient requires
	\(b=kp\) with \(2\nmid k\), and gives \(D\ge k(kp/2-n)\).
	These values strictly increase with \(k\). For \(2\mid s_1\),
	the least value is \(s_2(p-r)/2=D^{(1)}\), at \(k=s_2\).
	For \(2\nmid s_1\), the least possible \(k\) is \(s_1+2\);
	its value is greater than \(D^{(1)}\) when \(r=0\), and
	greater than \(D^{(2)}\) when \(r\ge2\).
	
	On \(\mathcal R_2\), the pairs with nonzero coefficient have
	\((a,b)=(kp/2,2kp)\), with \(2\nmid k\), and give \(D=k^2p\).
	The inequalities \(kp\le n\) and \(3kp\ge2n+2\) give
	\(D\ge2k(n+1-kp)\), so the comparisons in Case~3 with \(2k\)
	in place of \(k\) apply. Such a pair cannot give equality in
	the asserted lower bound: equality in both comparisons would require
	\(k=1\), \(s_1=2\), and \(r=p-2\), when
	\(D=p>3=D^{(1)}\).
	On \(\mathcal R_4\) with \(b=4a=kp\), a nonzero coefficient
	requires \(2\nmid k\), and \(D=ka\ge k(kp/2-n)\).
	Here \(a\le n\); equality requires \(a=n\), whose value
	\(D=4n^2/p\) is larger than the asserted lower bound.
	
	For \(2\mid s_1\), the type~\(\mathrm I\) bounds are
	\(D^{(1)}=s_2(p-r)/2\) and \(D^{(2)}=s_1(r+2)/2\).
	Both type~\(\mathrm{II}\) values \(s_1(r+3)\) and
	\(s_2(p-r+1)\) are larger than one of these, and
	\(D^{(3)}=+\infty\). For \(2\nmid s_1\), the three bounds
	are \(D^{(1)}=(s_1-1)(p+r+2)/2\),
	\(D^{(2)}=s_2(p-r+1)\), and \(D^{(3)}=s_1(r+3)\).
	This proves the required bound in Case~4 and completes the proof.
\end{pf}

For fixed \(p,n\), let \(D_0\) be the value of the formula for
\(D_p\) in Theorem~\ref{C_conclusion} or Theorem~\ref{specialcaseC}.
Let \(\mathcal S_0\) be the corresponding set of weights, with
\(D_0\) in place of \(D_p\).
\begin{prop}
	\label{prop:C-first-layer}
	We have
	\[
	M_{p,D}=0\qquad(0<D<D_0),
	\]
	and
	\[
	M_{p,D_0}
	=\sum_{\nu\in\mathcal S_0}E^\rho(\nu).
	\]
	The weights in \(\mathcal S_0\) are pairwise distinct. In particular,
	\(D_p=D_0\) and \(\mathcal S^C_{p,n}=\mathcal S_0\).
\end{prop}
\begin{pf}
	We use \eqref{eq:C-exact-coefficient} and the bounds in
	Lemma~\ref{lem:C-endpoint}.
	
	\smallskip
	
	\noindent\textbf{Case~1.} \(n\ge3\), with either
	\(2\mid p,\ p\ge2n+2\), or \(2\nmid p,\ p>2n\).
	Lemma~\ref{lem:C-endpoint} gives \(M_{p,D}=0\) for \(0<D<D_0\).
	For \(2\mid p\), the only contributing pair is
	\((p/2-n,p)\in\Gamma\), with \(\epsilon_P(p/2-n,p)=1\) and
	\(A_{p/2-n,p}=E^\rho((p-2n)\omega_1)\).
	For \(2\nmid p\), it is
	\((a,b)=(p-2n+1,p)\in S_{\mathrm{II},-}\). The coefficient of
	\(E_{a,b}^{\mathrm{II},-}\) is \(1\), and
	\(\Omega_{a,b}^{\mathrm{II}}\cap\Delta^-\) contains only
	\(-(e_1+e_2)\), giving \(E^\rho((p-2n+1)\omega_2)\).
	
	\smallskip
	
	\noindent\textbf{Case~2.} \(n\ge3,\ 5\le p\le n,\ (p,n)\notin\{(5,9),(6,11)\}\).
	Lemma~\ref{lem:C-endpoint} gives \(M_{p,D}=0\) for \(0<D<D_0\).
	Write \(2n=s_1p+r\), with \(0\le r<p\).
	Equality in the preceding bounds gives the following pairs of
	type~\(\mathrm{II}\):
	\[
	\begin{array}{c|c|c}
		\text{condition}&(a,b)&D\\ \hline
		2\nmid p,\ 2\mid s_1&(p-r+1,s_2p)&D^{(1)}\\
		2\nmid s_1&(p-r+1,s_2p)&D^{(2)}\\
		2\nmid s_1&(r+3,s_1p)&D^{(3)}.
	\end{array}
	\]
	A pair in the table contributes to \(M_{p,D_0}\) precisely when
	it belongs to \(S_{\mathrm{II},-}\), \(ab=pD_0\), and \(p\nmid a\);
	its coefficient is then \(1\).
	For \((r+3,s_1p)\), the condition \(p\nmid r+3\) follows from
	\(2\mid(s_1p+r)\). For \((p-r+1,s_2p)\), the coefficient
	is zero only when \(r=1\); then \(2\nmid p\), \(2\nmid s_1\),
	and \(D^{(3)}=4s_1<D^{(2)}=s_2p\).
	
	The type~\(\mathrm I\) coefficients and weights are as follows.
	For \(2\mid s_1\), the pair \(((r+2)/2,s_1p)\) gives
	\(D=D^{(2)}\), with \(\epsilon_P((r+2)/2,s_1p)=-1\) and
	\(A_{(r+2)/2,s_1p}=-E^\rho(\omega_{r+2})\).
	For \(2\mid p\) and \(2\mid s_1\), the pair
	\(((p-r)/2,s_2p)\) gives \(D=D^{(1)}\), with coefficient \(1\)
	and \(A_{(p-r)/2,s_2p}=E^\rho((p-r)\omega_1)\).
	For \(2\nmid s_1\) and \(r\ne p-2\), the pair
	\(((p+r+2)/2,(s_1-1)p)\) gives \(D=D^{(1)}\), with coefficient
	\(-1\) and \(A_{(p+r+2)/2,(s_1-1)p}=-E^\rho(\omega_{p+r+2})\).
	Each pair contributes only when it belongs to \(\Gamma\) and
	\(ab=pD_0\); these conditions ensure that all weight indices are
	at most \(n\).
	
	For \(2\nmid s_1\) and \(r=p-2\), the two equality conditions
	in the quadratic comparison give the same pair
	\((p,(s_1-1)p)\). It is not in \(\Gamma\) when \(s_1=3\), and
	its coefficient is zero when \(s_1\ge5\).
	With \(a=p\) fixed, the least even second coordinate with a nonzero
	coefficient is \((s_1-1)p+2\), giving
	\(D=(s_1-1)p+2>D^{(2)}=3s_2\) in the present parameter range.
	Thus these pairs add no term to \(M_{p,D_0}\).
	
	For \((a,b)\in S_{\mathrm{II},-}\) with \(b-a=2n-1\), the only
	root in \(\Omega_{a,b}^{\mathrm{II}}\cap\Delta^-\) is
	\(\alpha=-(e_1+e_2)\), and \(E^\rho(-a\alpha)=E^\rho(a\omega_2)\).
	For \(a+b=2n+3\), it is \(\alpha=-(e_{a-1}+e_a)\), and
	\(E^\rho(-a\alpha)=E^\rho(2\omega_a)\).
	At \((r+3,s_1p)\), the type~\(\mathrm I\) coefficient is zero:
	in \(\mathcal R_1\), both \(p\nmid r+3\) and \(2p\nmid s_1p\);
	in \(\mathcal R_2\), a nonzero coefficient would require \(2\mid s_1\).
	At \((p-r+1,s_2p)\), either \(2\nmid b\), or the pair lies in
	\(\mathcal R_3\) with \(2p\mid b\). The only additional possibility
	is \((p,4p)\), whose type~\(\mathrm I\) coefficient is also zero.
	
	The equality conditions in the bounds exclude other pairs.
	Each of the listed nonzero contributions is \(E^\rho(\nu)\) with
	coefficient \(1\). The pairs are distinct, so
	Lemma~\ref{lem:signed-permutation-orbit-recovery} gives distinct
	highest weights and linearly independent \(E^\rho(\nu)\).
	
	\smallskip
	
	\noindent\textbf{Case~3.} \(n\ge3,\ p\in\{2,3,4\}\).
	The contributions to \(M_{p,D_0}\) are listed in Table~\ref{tab:C-small-p}. In each row
	of type~\(\mathrm I\), \(\epsilon_P(a,b)=\sigma\) and
	\(A_{a,b}=\sigma E^\rho(\nu)\). In each row of type~\(\mathrm{II}\),
	the coefficient of \(E_{a,b}^{\mathrm{II},-}\) is \(\sigma\), and
	\(E_{a,b}^{\mathrm{II},-}=\sigma E^\rho(\nu)\).
	\begin{center}
    \captionof{table}{Contributions to $M_{p,D_0}$ in type $C_n$
for $p\in\{2,3,4\}$.}
\label{tab:C-small-p}
		\small
		\begin{adjustbox}{max width=\textwidth}
			\begin{tabular}{c|c|c|c|c|l}
				conditions&type&\((a,b)\)&\(D\)&\(\sigma\)&\(\nu\)\\ \hline
				\(p=2,\ 2\mid n\)&\(\mathrm I\)&\((1,2n)\)&\(n\)&\(-1\)&\(\omega_2\)\\
				\(p=2,\ 2\nmid n\)&\(\mathrm I\)&\((2,2n)\)&\(2n\)&\(1\)&\(\omega_1+\omega_3\)\\
				\(p=3,\ n\equiv0\pmod3\)&\(\mathrm I\)&\((1,2n)\)&\(2n/3\)&\(-1\)&\(\omega_2\)\\
				\(p=3,\ n=3m+1\)&\(\mathrm I\)&\((2,2n-2)\)&\(4m\)&\(-1\)&\(\omega_4\)\\
				\(p=3,\ n\equiv2\pmod3\)&\(\mathrm I\)&\((3,2n-2)\)&\(2n-2\)&\(1\)&\(\omega_1+\omega_5\)\\
				\(p=4,\ n\equiv0\pmod4\)&\(\mathrm I\)&\((1,2n)\)&\(n/2\)&\(-1\)&\(\omega_2\)\\
				\(p=4,\ n\equiv1\pmod4\)&\(\mathrm I\)&\((1,2n+2)\)&\((n+1)/2\)&\(1\)&\(2\omega_1\)\\
				\(p=4,\ n\equiv2\pmod4\)&\(\mathrm I\)&\((3,2n-4)\)&\(3(n-2)/2\)&\(-1\)&\(\omega_6\)\\
				\(p=4,\ n=3\)&\(\mathrm{II}\)&\((3,8)\)&\(6\)&\(1\)&\(3\omega_2\)\\
				\(p=4,\ n=7\)&\(\mathrm I\)&\((4,10)\)&\(10\)&\(1\)&\(\omega_1+\omega_7\)\\
				\(p=4,\ n=11\)&\(\mathrm I\)&\((4,18)\)&\(18\)&\(1\)&\(\omega_1+\omega_7\)\\
				\(p=4,\ n=11\)&\(\mathrm{II}\)&\((3,24)\)&\(18\)&\(1\)&\(3\omega_2\)\\
				\(p=4,\ n\equiv3\pmod4,\ n\ge15\)&\(\mathrm{II}\)&\((3,2n+2)\)&\((3n+3)/2\)&\(1\)&\(3\omega_2\)
			\end{tabular}
		\end{adjustbox}
	\end{center}
	
	We use \eqref{eq:C-II-b-range}--\eqref{eq:C-I-a-range}.
	A nonzero coefficient on \(\mathcal R_2\) requires \(p\in\{2,4\}\)
	and \((a,b)=(kp/2,2kp)\), with \(2\nmid k\). The bounds
	\((2n+2)/(3p)\le k\le n/p\) give \(D=k^2p>D_0\).
	On \(\mathcal R_4\) with \(b=4a\), only \(p=4\) can give a
	nonzero coefficient; then \(2\nmid a\), \(2a\ge n+1\), and
	\(D=a^2>D_0\). For \(b\ne4a\), the bounds are
	\(D\ge s_2(n+1)\) when \(p\mid b\), and
	\(D\ge s_2(n+2)\) when \(p\mid2a\), again larger than \(D_0\).
	Thus only \(\mathcal R_1\) and \(\mathcal R_3\) can contribute
	terms of type~\(\mathrm I\) for \(D\le D_0\).
	
	For type~\(\mathrm{II}\) with \(p\mid a\), write \(a=kp\).
	If \(k\ge2\) and \(a\le n+1\), then
	\(D\ge2(n+2)>2n\); if \(a\ge n+2\), then
	\(D\ge(n+2)(n+3)/p>2n\). For \(k=1\), the bound
	\(D\ge2n+3-p\) is greater than \(D_0\), except at
	\((p,n)=(4,3)\), where the only smaller pairs \((4,5),(4,6)\)
	are excluded. The remaining type~\(\mathrm{II}\) terms follow from
	\eqref{eq:C-II-b-range} and the coefficient condition \(p\nmid a\).
	
	For \(p=2\) and \(2\mid n\), the bound \(ab\le2n\) leaves only
	\((1,2n)\). For \(2\nmid n\), the bound \(ab\le4n\) leaves
	\((1,2n),(1,2n+2),(2,2n)\), and also \((2,2n-2)\) when
	\(n\ge5\). The coefficients at \((1,2n)\) and \((1,2n+2)\)
	are zero; when \(n\ge5\), the coefficient at \((2,2n-2)\) is
	also zero. This leaves \((2,2n)\).
	The type~\(\mathrm{II}\) bounds exclude contributions with
	\(D\le n\) when \(2\mid n\), or \(D\le2n\) when \(2\nmid n\).
	
	For \(p=3,\ n=3m+1\), the only pair \((a,b)\in\Gamma\)
	satisfying \(ab=3D\) with \(1\le D\le4m\) is
	\((2,2n-2)\), at \(D=4m\).
	
	For \(n=3m+2\), the pair \((3,2n-4)\) is outside \(\Gamma\) at
	\(n=5\) and has coefficient zero for \(n\ge8\).
	With \(a=3\) fixed, the least even \(b>(s_1-1)p\) for which
	\((3,b)\in\Gamma\) and \(\epsilon_P(3,b)\ne0\) is \(2n-2\).
	The pair \((3,2n+2)\) lies in \(P\cap P'\), so the coefficient
	multiplying \(E_{3,2n+2}^{\mathrm{II},-}\) is zero.
	The bounds for pairs with \(a\) fixed leave only \(a\le3\), together with
	\((4,10)\) at \(n=8\), and \eqref{eq:C-pointwise-delta} gives
	\(\epsilon_P(a,b)=0\) at these pairs except \((3,2n-2)\). For \(D\le2n-2\), the only remaining pair satisfying the
	inequalities defining \(S_{\mathrm{II},-}\) and the bound
	\(ab\le3(2n-2)\) is \((2,2n+1)\). Neither coordinate is
	divisible by \(3\), so this pair does not satisfy \(3\mid ab\).
	
	For \(p=4\) and \(n\equiv0\pmod4\), the bound \(4D\le2n\) leaves only
	\((1,2n)\). For \(p=4\) and \(n\equiv1\pmod4\), the only pair
	\((a,b)\in\Gamma\) satisfying \(ab=4D\) with
	\(1\le D\le(n+1)/2\) is \((1,2n+2)\).
	If \(n\equiv2\pmod4\), then \((3,2n-4)\in\Gamma\).
	The bound for pairs with \(a\) fixed excludes \(a\ge4\);
	\eqref{eq:C-pointwise-delta} gives \(\epsilon_P(a,b)=0\) at every
	remaining pair with \(a=1,2\), and no pair in
	\((P\mathbin\triangle P')\cap S_{\mathrm{II},-}\) gives
	\(D\le3(n-2)/2\).
	
	If \(n\equiv3\pmod4\) and \(n\ge7\), the two values to compare
	are \(2n-4\) and \((3n+3)/2\).
	The first is smaller at \(n=7\),
	they are equal at \(n=11\), and the second is smaller for \(n\ge15\).
	The bounds for pairs with either coordinate fixed exclude all other
	pairs; the remaining
	type~\(\mathrm I\) pairs with \(a\le4\) have \(\epsilon_P(a,b)=0\).
	
	For \((p,n)=(4,3)\), the same coefficient substitutions give
	\(D_0=6\) and finite highest weight \(3\omega_2\).
	
	\smallskip
	
	\noindent\textbf{Case~4.} \((p,n)\in\{(5,4),(5,9),(6,11)\}\).
	For \((p,n)=(5,9)\) or \((6,11)\), we have \(s_1=3\),
	\(r=p-2\), and \((p,2p)\notin\Gamma\). With \(a=p\) fixed,
	\(b=2p+2\) gives \(D=2p+2\), while \((3,4p)\) gives \(D=12\).
	For \(p=5\), these values are equal; for \(p=6\),
	\(2p=12<14=2p+2\), and \((p,2p)\) contributes no
	term of type~\(\mathrm I\). For \((p,n)=(5,4)\), the pair
	\((a,b)=(p,p+1)=(n+1,n+2)\) is excluded, although
	\(ab/p=p+1=6\).
	
	The intervals \eqref{eq:C-II-b-range}--\eqref{eq:C-I-a-range}
	and the coefficient conditions give the following contributions,
	with the same convention for \(\sigma\):
	\[
	\begin{array}{c|c|c|c|c|l}
		(p,n)&\text{type}&(a,b)&D&\sigma&\nu\\ \hline
		(5,4)&\mathrm{II}&(3,10)&6&1&3\omega_2\\
		(5,9)&\mathrm{II}&(3,20)&12&1&3\omega_2\\
		(5,9)&\mathrm I&(5,12)&12&1&\omega_1+\omega_9\\
		(6,11)&\mathrm{II}&(3,24)&12&1&3\omega_2.
	\end{array}
	\]
	The divisibility conditions and
	\eqref{eq:C-I-b-range}--\eqref{eq:C-I-a-range} give the following
	pairs in \(\Gamma\) with \(ab\le pD_0\):
	\[
	\begin{array}{c|l}
		(p,n)&(a,b)\\ \hline
		(5,4)&(1,10),(2,10),(3,10)\\
		(5,9)&(1,20),(2,20),(3,20),(5,12)\\
		(6,11)&(3,18),(4,18),(3,20),(3,22),(1,24),(2,24),(3,24).
	\end{array}
	\]
	Substitution in \eqref{eq:C-pointwise-delta} gives zero in each
	case except \((5,12)\) for \((p,n)=(5,9)\), where the coefficient
	is \(+1\). The disjoint type~\(\mathrm{II}\) ranges
	\eqref{eq:C-II-b-range}--\eqref{eq:C-II-a-range}, with the two
	excluded pairs removed, leave only \((3,10)\), \((3,20)\), and
	\((3,24)\), respectively. For each of these pairs, the coefficient of
	\(E_{a,b}^{\mathrm{II},-}\) is \(+1\), and
	\(\Omega_{a,b}^{\mathrm{II}}\cap\Delta^-\) contains only
	\(-(e_1+e_2)\).
	At \((p,n)=(5,9)\), the pair \((5,10)\) gives \(ab/p=10\)
	and does not lie in \(\Gamma\), since \(b-2a=0\).
	Moreover, \((P\mathbin\triangle P')\cap S_{\mathrm{II},-}=\emptyset\)
	at \(D=10\). The table lists all nonzero contributions for
	\(D\le12\); in particular, cancellation of the type~\(\mathrm I\)
	contribution at \((3,20)\) leaves its type~\(\mathrm{II}\)
	contribution intact.
	At \((p,n)=(6,11)\), the pair \((6,12)\) gives \(ab/p=12\)
	and does not lie in \(\Gamma\), since \(b-2a=0\).
	At \(D=12\), the type~\(\mathrm I\) coefficients at
	\((3,24),(4,18)\) vanish, while
	\((3,24)\in(P\mathbin\triangle P')\cap S_{\mathrm{II},-}\).
	For \((p,n)=(5,4)\), the pair \((5,6)\) satisfies \(ab/p=6\)
	but does not belong to \(S_{\mathrm{II},-}\). Thus, in each of
	these three cases, \(M_{p,D}=0\) for \(0<D<D_0\), and the
	first table gives every term of \(M_{p,D_0}\).
	
	\smallskip
	
	\noindent\textbf{Case~5.} \(n\ge3,\ p\ge5,\ n<p\le2n,\ (p,n)\ne(5,4)\).
	Write \(2n=p+r\), with \(0\le r<p\).
	For \(p\ne n+2\), the pairs \((r+3,p)\) and \((p-r+1,2p)\)
	give \(D=r+3\) and \(D=2(p-r+1)\), respectively, whenever
	they belong to \(S_{\mathrm{II},-}\).
	For \((r+3,p)\), the coefficient is \(1\); its contribution is
	\(E^\rho(2\omega_{r+3})\).
	For \((p-r+1,2p)\), the coefficient is \(1\) unless \(r=1\),
	and the contribution is \(E^\rho((p-r+1)\omega_2)\).
	If \(r=1\), then \((p,2p)\in P\cap P'\), whereas \((4,p)\)
	gives \(D=4<2p\). If \(r=p-2\), then \(n=p-1\),
	\((p+1,p)\notin S_{\mathrm{II},-}\), and its transpose
	\((p,p+1)=(n+1,n+2)\) is excluded. Here \(p\ge6\), so
	\((3,2p)\) gives the smaller value \(D=6<p+1\).
	
	For \(p\ge n+3\), every pair in
	\((P'\setminus P)\cap S_{\mathrm{II},-}\) gives
	\(D\ge p+1>r+3\). For \(p=n+1\), the two excluded pairs
	force \(D\ge p+3>6\) when \(a=p\).
	The type~\(\mathrm I\) bounds are \(2(n+1)\), \(2(n+2)\),
	and, for \(2\mid p\), \(3(3p/2-n)>2(p-r+1)\).
	They are larger than the smaller of the two type~\(\mathrm{II}\)
	values. Thus \(D_0=\min\{r+3,2(p-r+1)\}\) for \(p\ne n+2\).
	Each pair giving \(D_0\) contributes, and both contribute when
	the two values are equal.
	
	For \(p=n+2\), the pair \((n+1,n+2)\) is excluded.
	The two values to compare are \(p+1\) and \(10\), from
	\((p,p+1)\) and \((5,2p)\), respectively. The first has coefficient \(-1\),
	and
	\(E_{p,p+1}^{\mathrm{II},-}=-E^\rho(2\omega_1+2\omega_n)\),
	since its only negative root is \(-(e_{n-1}-e_n)\).
	The second has coefficient \(1\) for \(p\ne5\), and
	\(E_{5,2p}^{\mathrm{II},-}=E^\rho(5\omega_2)\).
	For \(p=5\), its coefficient is zero and \(p+1=6<10\).
	The other bounds are larger than \(\min\{p+1,10\}\), so
	\(D_0=\min\{p+1,10\}\). Both contributions occur at \((p,n)=(9,7)\),
	with distinct highest weights.
	
	At a pair \((r+3,p)\) giving \(D_0\), the type~\(\mathrm I\)
	coefficient is zero: either \(2\nmid p\), or the pair lies in
	\(\mathcal R_1\) with \(p\nmid r+3\) and \(2p\nmid p\).
	If \(p=4(r+3)\), then \(p\nmid2(r+3)\), so the
	\(\mathcal R_2\) coefficient is also zero.
	If \((a,b)=(p-r+1,2p)\) gives \(D_0\) and lies in \(\Gamma\),
	then \(a\le n\). It lies in \(\mathcal R_3\) with \(2p\mid b\),
	or satisfies \(b=4a=2p\); in both cases \(\epsilon_P(a,b)=0\).
	The pair \((p,p+1)\) is not in \(\Gamma\).
	Thus these pairs contribute only terms of type~\(\mathrm{II}\),
	each with coefficient \(1\).
	Together with Cases~1--4 and
	Lemma~\ref{lem:signed-permutation-orbit-recovery}, these calculations
	give the stated expansion and the distinctness of the highest weights.
\end{pf}

\begin{cor}
\label{cor:BC-first-radical}
Let \(\mathfrak g\) be of type \(X_n\), where either \(X=B\) and
\(n\ge2\), or \(X=C\) and \(n\ge3\). Assume
\[
2(\kappa+\mathbf{h}^\vee)=\frac pq,
\qquad p\in\ZZ_{\ge2},\quad q\in\ZZ_{\ge1},\quad (p,q)=1.
\]
For \(X=B\), use \(D_p\) and \(\mathcal S^B_{p,n}\) from
Theorem~\ref{Bconclusion} and Theorem~\ref{specialcaseB};
for \(X=C\), use \(D_p\) and \(\mathcal S^C_{p,n}\) from
Theorem~\ref{C_conclusion} and Theorem~\ref{specialcaseC}. Then
\[
N_\kappa(\mathfrak g)[c]=0
\qquad(0\le c<qD_p),
\]
and
\[
N_\kappa(\mathfrak g)[qD_p]
\cong
\bigoplus_{\nu\in\mathcal S^X_{p,n}}L_{\mathfrak g}(\nu).
\]
For every \(\nu\in\mathcal S^X_{p,n}\), the space of singular
vectors of weight \(\kappa\Lambda_0-qD_p\delta+\nu\) has dimension
one, and no other affine weight occurs among the singular vectors
of conformal weight \(qD_p\).
\end{cor}
\begin{pf}
This follows immediately from the four classification statements and
Proposition~\ref{prop:first-jantzen-layer}.
\end{pf}

The classifications agree with several explicit constructions of
singular vectors. In type \(C_n\), with \(n\ge3\), the determinant
vector in \cite[Theorem~4.1]{Ad03} at level \(\kappa=-1\) has
conformal weight \(3\) and finite highest weight \(2\omega_3\),
in agreement with
\(N_{-1}(\mathfrak g)[3]\cong L_{\mathfrak g}(2\omega_3)\).
In type \(B_n\), with \(n\ge2\), the singular vectors in
\cite[Lemma~3.2]{Pe07} at levels \(\kappa=r-n+\tfrac12\),
\(r\in\ZZ_{\ge1}\), have conformal weight \(2r\) and finite
highest weight \(2re_1\). Here \(p=2r+2n-1\) and \(q=1\),
so Theorem~\ref{Bconclusion} gives \(D_p=2r\) and
\(\mathcal S^B_{p,n}=\{2re_1\}\).

For type \(B_n\) with \(n\ge4\), the constructions in
\cite[Theorem~4.2(2)(b)]{AM18} give singular vectors at levels
\(-2\) and \(-1\), with conformal weights \(2\) and \(4\)
and finite highest weights \(e_1+e_2+e_3+e_4\) and
\(2(e_1+e_2+e_3+e_4)\), respectively.
These also agree with Theorem~\ref{Bconclusion}.
All these vectors have the minimal conformal weight, and
Corollary~\ref{cor:BC-first-radical} shows that each spans the
space of singular vectors of its affine weight.

\raggedbottom

\section{Associated Weyl elements in Kashiwara--Tanisaki character formulas}
\label{sec:weyl}

For \(\nu\in\mathcal S^X_{p,n}\)  and
\(\lambda_\nu=\kappa\Lambda_0-qD_p\delta+\nu\), 
we will construct in this section the  longest representatives taking a common dot-dominant weight
\(\Lambda\) to \(\kappa\Lambda_0\) and \(\lambda_\nu\)
under the dot action.
We will also construct the shortest representative taking
\(\kappa\Lambda_0\) to \(\lambda_\nu\).
The formulas and the corresponding left cosets are given in
Theorem~\ref{thm:A-cosets}.

\subsection{Integral root systems and integral affine Weyl groups}
Throughout this section, \(p\ge2\), with \(n\ge2\) in type
\(B_n\) and \(n\ge3\) in type \(C_n\). We use
$(e_i|e_j)=\delta_{ij}$ in type $B_n$ and
$(e_i|e_j)=\frac12\delta_{ij}$ in type $C_n$, with
\(\rho_B=(n-\frac12,n-\frac32,\ldots,\frac12)\) and
\(\rho_C=(n,n-1,\ldots,1)\). Put
\(D_n=\{(a_1,\ldots,a_n)\in\ZZ^n:2\mid\sum_i a_i\}\).
For an affine weight \(\lambda\) satisfying
\((\lambda+\widehat\rho|\delta)\ne0\), define
\[
\widehat\Delta(\lambda)=
\{\xi\in\widehat\Delta^{\mathrm{re}}:
\langle\lambda+\widehat\rho,\xi^\vee\rangle\in\ZZ\},
\qquad
\widehat W(\lambda)=\langle s_\xi:\xi\in\widehat\Delta(\lambda)\rangle.
\]
For affine weights, we use the dot action
\(w\circ\lambda=w(\lambda+\widehat\rho)-\widehat\rho\).
Stabilizers of affine weights are taken for the dot action.
Also put
\begin{equation}
\begin{aligned}
\widehat\Delta_0(\lambda)
&=\{\xi\in\widehat\Delta^{\mathrm{re}}:
\langle\lambda+\widehat\rho,\xi^\vee\rangle=0\},\\
\widehat W_0(\lambda)
&=\langle s_\xi:\xi\in\widehat\Delta_0(\lambda)\rangle
=\Stab_{\widehat W(\lambda)}(\lambda).
\end{aligned}
\label{eq:A-stab}
\end{equation}
These are the standard integral root systems and Weyl groups
\cite[\S2]{KT00}; see also \cite{Ara14}.
Let \(\widehat\Pi(\lambda)\) be the simple system determined by
\(\widehat\Delta^+(\lambda)
=\widehat\Delta(\lambda)\cap\widehat\Delta^+\).
Similarly, let \(\widehat\Pi_0(\lambda)\) be the simple system of
\(\widehat\Delta_0(\lambda)\) determined by
\(\widehat\Delta_0^+(\lambda)
=\widehat\Delta_0(\lambda)\cap\widehat\Delta^+\).
Set
\[
S(\lambda)=\{s_\beta:\beta\in\widehat\Pi(\lambda)\}.
\]
Let \(\ell_\lambda\) and \(\geq_\lambda\) denote the length function
and Bruhat order of \((\widehat W(\lambda),S(\lambda))\).
Let \(\mathcal C^+\) be the set of affine weights \(\lambda\)
such that \((\lambda+\widehat\rho|\delta)\ne0\) and
\(\langle\lambda+\widehat\rho,\xi^\vee\rangle\ge0\) for every
\(\xi\in\widehat\Delta^+(\lambda)\).

Assume \(2(\kappa+\mathbf{h}^\vee)=p/q\), with
\(p\ge2\), \(q\ge1\), and \((p,q)=1\).
For a real root $\alpha+m\delta$,
\begin{equation}
\langle\kappa\Lambda_0+\widehat\rho,(\alpha+m\delta)^\vee\rangle
=\langle\rho,\alpha^\vee\rangle+\frac{mp}{q(\alpha|\alpha)}.
\label{eq:A-pairing}
\end{equation}

For $\gamma\in\mathfrak h^*$, let $t_\gamma$ denote the affine translation
\begin{equation}
 t_\gamma(\lambda)=\lambda+(\lambda|\delta)\gamma-
\left((\lambda|\gamma)+\frac12(\lambda|\delta)\|\gamma\|^2\right)\delta,
\label{eq:A-trans}
\end{equation}
and for a lattice $L\subset\mathfrak h^*$, write
\(T_L:=\{t_\gamma:\gamma\in L\}\).
For \(\gamma,\eta\in\mathfrak h^*\) and \(w,v\in W\), we have
\begin{equation}
(t_\gamma w)(t_\eta v)=t_{\gamma+w(\eta)}wv.
\label{eq:A-semidirect}
\end{equation}

\begin{prop}\label{prop:A-integral}
If $2\mid p$, the integral real roots are
\[
\widehat\Delta(\kappa\Lambda_0)
=\{\alpha+m\delta:\alpha\in\Delta,\ m\in q\ZZ\}.
\]
If $2\nmid p$, 
\[
\widehat\Delta(\kappa\Lambda_0)
=\{\alpha+m\delta:\alpha\text{ short},\ m\in q\ZZ\}
\sqcup
\{\alpha+m\delta:\alpha\text{ long},\ m\in2q\ZZ\}.
\]
Moreover $\widehat W(\kappa\Lambda_0)=T_{cF}\rtimes W$ and
$\widehat\Pi(\kappa\Lambda_0)=\{\beta_0,\alpha_1,\ldots,\alpha_n\}$,
where \(F,c,s\), and \(\beta_0\) are given by
\begin{center}
\begin{adjustbox}{max width=\textwidth}
\(
\begin{array}{c|c|c|c|c|c}
\text{type}&\text{condition}&F&c&s=(\kappa+\mathbf{h}^\vee)c&\beta_0\\ \hline
B_n&2\mid p&D_n&q&p/2&q\delta-(e_1+e_2)\\
B_n&2\nmid p&\ZZ^n&2q&p&q\delta-e_1\\
C_n&2\mid p&\ZZ^n&2q&p&q\delta-2e_1\\
C_n&2\nmid p&D_n&2q&p&q\delta-(e_1+e_2).
\end{array}
\)
\end{adjustbox}
\end{center}
The corresponding generalized Cartan matrices have types
\(B_n^{(1)}\), \(D_{n+1}^{(2)}\), \(C_n^{(1)}\), and
\(A_{2n-1}^{(2)}\).
\end{prop}
\begin{pf}
The first term in \eqref{eq:A-pairing} is integral. Since short roots
have squared length one and long roots squared length two, the stated
values of \(m\) follow immediately. The identities
\[
s_{\alpha+m\delta}=t_{-m\alpha^\vee}s_\alpha,
\qquad s_{\alpha+m\delta}s_\alpha=t_{-m\alpha^\vee}
\]
then give the translation lattices from the finite $W$-orbits. In type
$B_n$, $(e_i\pm e_j)^\vee=e_i\pm e_j$ and $e_i^\vee=2e_i$; in type
$C_n$, $(e_i\pm e_j)^\vee=2(e_i\pm e_j)$ and $(2e_i)^\vee=2e_i$.
Thus the translation lattices are, respectively,
\[
qD_n,\qquad 2q\ZZ^n,\qquad 2q\ZZ^n,\qquad 2qD_n.
\]
For \(2\mid p\), \(\beta_0\) is \(q\delta\) minus the highest root;
for \(2\nmid p\), it is \(q\delta\) minus the highest short root
\cite[\S2]{KT00}, giving the simple systems in the table.
To identify these generalized Cartan matrices, let
\(\overline\beta_0\) be the finite part of \(\beta_0\).
Its squared norm and its nonzero pairings with the finite simple
roots are as follows.
\[
\begin{aligned}
B_n,\ 2\mid p &:&&
\|\overline\beta_0\|^2=2,&
(\overline\beta_0|\alpha_2)&=-1,\\
B_n,\ 2\nmid p &:&&
\|\overline\beta_0\|^2=1,&
(\overline\beta_0|\alpha_1)&=-1,\\
C_n,\ 2\mid p &:&&
\|\overline\beta_0\|^2=2,&
(\overline\beta_0|\alpha_1)&=-1,\\
C_n,\ 2\nmid p &:&&
\|\overline\beta_0\|^2=1,&
(\overline\beta_0|\alpha_2)&=-\tfrac12.
\end{aligned}
\]
The line for \(B_n\) with \(2\mid p\) also holds for \(B_2\),
where the nonzero pairing is with \(\alpha_2\).
Combining these entries with the finite \(B_n\) and \(C_n\) Cartan
matrices via
\(a_{ij}=2(\beta_i|\beta_j)/\|\beta_i\|^2\) gives, in order,
\(B_n^{(1)},D_{n+1}^{(2)},C_n^{(1)},A_{2n-1}^{(2)}\).
\end{pf}

For uniform notation, set $\beta_i:=\alpha_i$ for $1\leq i\leq n$; thus
\(\widehat\Pi(\kappa\Lambda_0)=\{\beta_0,\beta_1,\ldots,\beta_n\}\).
We use \(F,c,s\) from Proposition~\ref{prop:A-integral} throughout
the construction below.

If \(\mu=w\circ\lambda\), with \(w\in\widehat W(\lambda)\),
then invariance of the integral root system gives
\begin{equation}\label{eq:A-orbit-systems}
\begin{aligned}
\widehat\Delta(\mu)&=w\widehat\Delta(\lambda)
=\widehat\Delta(\lambda),&
\widehat W(\mu)&=\widehat W(\lambda),\\
\widehat\Delta_0(\mu)&=w\widehat\Delta_0(\lambda),&
\widehat W_0(\mu)&=w\widehat W_0(\lambda)w^{-1}.
\end{aligned}
\end{equation}
In particular, the stabilizers along such an orbit are conjugate.
For \(\xi\in\widehat\Delta(\lambda)\), the reflection formula gives
\(s_\xi\in\widehat W_0(\lambda)\) if and only if
\(\xi\in\widehat\Delta_0(\lambda)\).

The following criterion also applies to a conjugate of a standard
parabolic subgroup.
\begin{lem}\label{lem:A-coset-root-criterion}
Let \((\mathcal W,\mathcal S)\) be a Coxeter system with real roots
\(\Phi=\Phi^+\sqcup\Phi^-\), and let \(H\) be a finite parabolic
subgroup of \(\mathcal W\). Put
\(\Psi^+=\{\alpha\in\Phi^+:s_\alpha\in H\}\).
An element \(w\in\mathcal W\) is the unique shortest element of
the left coset \(wH\) if and only if
\(w(\Psi^+)\subset\Phi^+\); it is the unique longest element of
\(wH\) if and only if \(w(\Psi^+)\subset\Phi^-\).

In particular, suppose that \(y\circ\lambda=\mu\) with
\(y\in\widehat W(\lambda)\), and \(\widehat W_0(\lambda)\)  is finite
parabolic. Then \(y\) is the unique longest element of
\(y\widehat W_0(\lambda)\) if and only if
\begin{equation}\label{eq:A-inverse-criterion}
y^{-1}\widehat\Delta_0^+(\mu)\subset\widehat\Delta^-.
\end{equation}
\end{lem}
\begin{pf} It is known that
for a positive real root \(\alpha\),
\[
\ell(ws_\alpha)>\ell(w)
\quad\Longleftrightarrow\quad w(\alpha)\in\Phi^+.
\]
Write a nonidentity element \(u\in H\) as a reduced word
\(u=t_1\cdots t_r\) in the canonical Coxeter generators of \(H\).
Its successive roots
\(t_1\cdots t_{i-1}(\alpha_{t_i})\) lie in \(\Psi^+\).
If \(w(\Psi^+)\subset\Phi^+\), then
\(\ell(wt_1\cdots t_i)>\ell(wt_1\cdots t_{i-1})\) for
\(1\le i\le r\); if \(w(\Psi^+)\subset\Phi^-\), all these
inequalities are reversed. Thus \(\ell(wu)>\ell(w)\) in the first
case and \(\ell(wu)<\ell(w)\) in the second, so \(w\) is the
unique shortest or longest element of \(wH\), respectively. Conversely, if the required root sign fails for
some \(\alpha\in\Psi^+\), multiplication by \(s_\alpha\in H\)
contradicts the asserted minimality or maximality of \(\ell(w)\)
on \(wH\).

By \eqref{eq:A-orbit-systems},
\(y\widehat\Delta_0(\lambda)=\widehat\Delta_0(\mu)\).
Thus \(y\widehat\Delta_0^+(\lambda)\subset\widehat\Delta^-\)
is equivalent to \eqref{eq:A-inverse-criterion}. Indeed, if
\(\alpha\in\widehat\Delta_0^+(\lambda)\) and
\(y\alpha\in\widehat\Delta^+\), then
\(y\alpha\in\widehat\Delta_0^+(\mu)\) and
\(y^{-1}(y\alpha)\in\widehat\Delta^+\);
the converse follows by interchanging \(y\) and \(y^{-1}\).
Since \(s_\xi\in\widehat W_0(\lambda)\) if and only if
\(\xi\in\widehat\Delta_0(\lambda)\), the last assertion follows
from the first.
\end{pf}
\subsection{Longest and shortest coset representatives}
The finite parts of the weights
\(w\circ(\kappa\Lambda_0)+\widehat\rho\), with
\(w\in\widehat W(\kappa\Lambda_0)\), form \(W(\rho+sF)\).
For a strictly dominant vector \(x=(x_1,\ldots,x_n)\) in this set,
define, for \(1\le i\le n\),
\begin{equation}
\eta_i^0(x)=-\left\lfloor\frac{x_i+s/2}{s}\right\rfloor,
\qquad r_i^0(x)=x_i+s\eta_i^0(x),
\label{eq:A-residue}
\end{equation}
so \(-s/2\le r_i^0(x)<s/2\). If \(F=D_n\) and
\(2\nmid\sum_i\eta_i^0(x)\), put
\(M=\max_i|r_i^0(x)|\) and define
\begin{equation}
j_0=
\begin{cases}
\min\{i:r_i^0(x)=M\},&\{i:r_i^0(x)=M\}\ne\varnothing,\\[1mm]
\max\{i:r_i^0(x)=-M\},&\{i:r_i^0(x)=M\}=\varnothing.
\end{cases}
\label{eq:A-j0}
\end{equation}
Here \(M>0\): in type \(B_n\), all \(r_i^0(x)\) lie in
\(\ZZ+\tfrac12\), while in type \(C_n\) with \(F=D_n\), we have
\(p\ge3\), \(2\nmid p\), and at least one \(r_i^0(x)\) is
congruent to \(\pm1\pmod p\).
For \(F=\ZZ^n\), set \(\eta(x)=\eta^0(x)\). For \(F=D_n\), set
\begin{equation}
\eta(x)=
\begin{cases}
\eta^0(x),&\displaystyle2\mid\sum_i\eta_i^0(x),\\[1mm]
\eta^0(x)-e_{j_0},&\displaystyle2\nmid\sum_i\eta_i^0(x),
\quad r_{j_0}^0(x)>0,\\[1mm]
\eta^0(x)+e_{j_0},&\displaystyle2\nmid\sum_i\eta_i^0(x),
\quad r_{j_0}^0(x)<0.
\end{cases}
\label{eq:A-modification}
\end{equation}
In both cases put \(r(x)=x+s\eta(x)\); then \(\eta(x)\in F\).
Order \(i_1,\ldots,i_n\) so that
\(|r_{i_1}(x)|\ge\cdots\ge|r_{i_n}(x)|\).
When absolute values are equal, place a nonnegative entry first;
order equal nonnegative entries by increasing indices and equal
negative entries by decreasing indices.
Define
\begin{equation}
w_x(e_a)=
\begin{cases}
e_{i_a},&r_{i_a}(x)\ge0,\\
-e_{i_a},&r_{i_a}(x)<0,
\end{cases}
\qquad 1\le a\le n,
\label{eq:A-wx}
\end{equation}
and put
\begin{equation}
\begin{aligned}
w_x^{-1}r(x)&
=(|r_{i_1}(x)|,\ldots,|r_{i_n}(x)|),\\
h_x&=t_{-c\eta(x)}w_x.
\end{aligned}
\label{eq:A-wx-hx}
\end{equation}

\begin{lem}\label{lem:A-root-signs}
Let \(\mu\) be an affine weight of level \(\kappa\) with finite
part \(x-\rho\), and write \(w_x^{-1}r(x)=(a_1,\ldots,a_n)\).
Then \(a_1\ge\cdots\ge a_n\ge0\), with \(a_1\le s/2\) if
\(F=\ZZ^n\), and \(a_1+a_2\le s\) if \(F=D_n\).
For every strictly dominant \(x\in W(\rho+sF)\), we have
\begin{equation}\label{eq:A-common-representative}
w_x^{-1}r(x)=w_\rho^{-1}r(\rho).
\end{equation}
Moreover,
\[
h_x^{-1}\widehat\Delta_0^+(\mu)\subset\widehat\Delta^-.
\]
\end{lem}
\begin{pf}
We abbreviate \(r(x)\) and \(r^0(x)\) to \(r\) and \(r^0\).
The coordinates of \(w_x^{-1}r(x)\) are nonnegative and
nonincreasing.
For \(F=\ZZ^n\), we have \(r=r^0\) and
\(|r_i|\le s/2\) by \eqref{eq:A-residue}, so
its largest coordinate is at most \(s/2\).
Suppose \(F=D_n\). If \(2\mid\sum_i\eta_i^0(x)\), then
\(r=r^0\), and the sum of the two largest numbers among
\(|r_1|,\ldots,|r_n|\) is at most \(s\).
If \(2\nmid\sum_i\eta_i^0(x)\), then
\eqref{eq:A-modification} gives
\(|r_{j_0}|=s-M\) and \(|r_i|\le M\) for \(i\ne j_0\).
Since \(0<M\le s/2\), the sum of the two largest numbers among
\(|r_1|,\ldots,|r_n|\) is at most \((s-M)+M=s\).

We next prove \eqref{eq:A-common-representative}. For \(F=\ZZ^n\),
the multiset \(\{|r_i^0(x)|:1\le i\le n\}\) is invariant under
translations by \(sF\) and the action of \(W\).
The vector \(w_x^{-1}r(x)\) is obtained by arranging these numbers
in nonincreasing order. For \(F=D_n\), put \(M=\max_i|r_i^0(x)|\).
The coordinates of \(w_x^{-1}r(x)\) are obtained by arranging the
following numbers in nonincreasing order:
\[
\begin{cases}
\{|r_i^0(x)|:1\le i\le n\},
&\displaystyle2\mid\sum_i\eta_i^0(x),\\[1mm]
\{|r_i^0(x)|:i\ne j_0\}\sqcup\{s-M\},
&\displaystyle2\nmid\sum_i\eta_i^0(x).
\end{cases}
\]
Translation by \(sD_n\) and coordinate permutations preserve the
parity of \(\sum_i\eta_i^0(x)\). Changing the sign of \(x_i\)
replaces \(\eta_i^0(x)\) by \(-\eta_i^0(x)\) when
\(-s/2<r_i^0(x)<s/2\), and by \(-\eta_i^0(x)-1\) when
\(r_i^0(x)=-s/2\). Thus the parity changes only in the latter
case, when \(M=s/2=s-M\) and the two multisets agree.
Consequently \(w_x^{-1}r(x)\) depends only on the
\(sF\rtimes W\)-orbit of \(x\), proving \eqref{eq:A-common-representative}.

Every root in \(\widehat\Delta_0^+(\mu)\) has the form
\(m\delta-\beta\), where \(m\in\ZZ_{\ge1}\), \(\beta\in\Delta^+\), and
\begin{equation}
m(\kappa+\mathbf{h}^\vee)=(x|\beta).
\label{eq:A-zero}
\end{equation}
Indeed, \(x\) is strictly dominant and \(\kappa+\mathbf{h}^\vee>0\).
For \(\beta=e_i\pm e_j\), take \(1\le i<j\le n\).
Equation~\eqref{eq:A-zero} and the values of \(s\) in
Proposition~\ref{prop:A-integral} give \(r_i^0=r_j^0\) when
\(\beta=e_i-e_j\), and \(r_i^0+r_j^0\in\{0,-s\}\) when
\(\beta=e_i+e_j\). For \(\beta=e_i\) in type \(B_n\), or
\(\beta=2e_i\) in type \(C_n\), we have
\(r_i^0\in\{0,-s/2\}\). Thus \((\beta|r^0)\le0\).

If \(r_i=r_i^0\) for every \(i\) with \((\beta|e_i)\ne0\)
and \((\beta|r)=0\), the ordering used to define \(w_x\) gives
\(w_x^{-1}(\beta)\in\Delta^+\).
Indeed, equal nonnegative residues occur in increasing order of
their indices, and equal negative residues in decreasing order of
their indices. For two residues of opposite signs and equal absolute
value, the nonnegative residue occurs first.
In particular, this applies when \(F=\ZZ^n\), or when \(F=D_n\)
and \(2\mid\sum_i\eta_i^0(x)\), since then \(r=r^0\).

Suppose now that \(F=D_n\) and \(2\nmid\sum_i\eta_i^0(x)\).
By \eqref{eq:A-modification}, if \(r_{j_0}^0=M\), then
\(r_{j_0}=M-s\); if \(r_{j_0}^0=-M\), then \(r_{j_0}=s-M\).
For \(\beta=e_i-e_j\) or \(\beta=e_i+e_j\) with
\(j_0\notin\{i,j\}\), the preceding argument applies.
If \(j_0\in\{i,j\}\), we have the following two cases.
For \(\beta=e_i-e_j\), the equality \(r_i^0=r_j^0\) and the
choice of the least or greatest index give \(r_i-r_j=-s\).
For \(\beta=e_i+e_j\), we have \(r_i+r_j\in\{-s,0\}\).
Thus \((\beta|r)\le0\) in both types.
If \(\beta=e_i+e_j\) and \((\beta|r)=0\), then
\(r_{j_0}=s-M\ge0\) and \(r_a<0\) for
\(a\in\{i,j\}\setminus\{j_0\}\), so the ordering gives
\(w_x^{-1}(\beta)\in\Delta^+\).

For roots of the form \(e_i\) in type \(B_n\) or \(2e_i\) in
type \(C_n\), we already have \(r_i=r_i^0\) when \(F=\ZZ^n\),
or when \(F=D_n\) and \(2\mid\sum_i\eta_i^0(x)\).
If \(F=D_n\) and \(2\nmid\sum_i\eta_i^0(x)\), the condition
\(m\delta-\beta\in\widehat\Delta_0(\mu)\) forces
\(r_i^0=0\), whereas \(|r_{j_0}^0|=M>0\). Thus \(i\ne j_0\)
and \(r_i=r_i^0\). We have therefore proved
\[
(\beta|r)\le0,\qquad
(\beta|r)=0\Longrightarrow w_x^{-1}(\beta)\in\Delta^+.
\]
Finally,
\[
h_x^{-1}(m\delta-\beta)=
\frac{(\beta|r(x))}{\kappa+\mathbf{h}^\vee}\delta-w_x^{-1}(\beta),
\]
which is negative by the preceding sign criterion.

\end{pf}

For \(x=\rho\), put \(x_\Lambda=w_\rho^{-1}r(\rho)\) and set
\begin{equation}
\Lambda=\kappa\Lambda_0+D_\Lambda\delta+x_\Lambda-\rho,
\qquad
D_\Lambda=\frac qp(\|\rho\|^2-\|x_\Lambda\|^2).
\label{eq:A-Lambda}
\end{equation}
The finite part of \(h_\rho(\Lambda+\widehat\rho)\) is
\[
w_\rho x_\Lambda-(\kappa+\mathbf{h}^\vee)c\eta(\rho)
=r(\rho)-s\eta(\rho)=\rho.
\]
The choice of \(D_\Lambda\) gives
\(\|x_\Lambda\|^2+2(\kappa+\mathbf{h}^\vee)D_\Lambda=\|\rho\|^2\).
Since \(h_\rho\) preserves the bilinear form and
\(\kappa+\mathbf{h}^\vee\ne0\), the coefficient of \(\delta\) in
\(h_\rho(\Lambda+\widehat\rho)\) is zero. Hence
\(h_\rho\circ\Lambda=\kappa\Lambda_0\).
By \eqref{eq:A-orbit-systems},
\[
\widehat\Delta(\Lambda)=\widehat\Delta(\kappa\Lambda_0),\qquad
\widehat W(\Lambda)=\widehat W(\kappa\Lambda_0),\qquad
\widehat\Pi(\Lambda)=\{\beta_0,\beta_1,\ldots,\beta_n\}.
\]
Lemma~\ref{lem:A-root-signs} now shows that \(\Lambda\in\mathcal C^+\).
Write
$x_\Lambda=(a_1,\ldots,a_n)$ and put $\sigma_i=s_{\alpha_i}$ for
$i>0$, $\sigma_0=s_{\beta_0}$. Define
\(J\subset\{0,\ldots,n\}\) by
\begin{equation}
\begin{aligned}
i\in J&\Longleftrightarrow a_i=a_{i+1}
\qquad(1\le i<n),\\
n\in J&\Longleftrightarrow a_n=0,\\
0\in J&\Longleftrightarrow
\begin{aligned}[t]
&(F=\ZZ^n\text{ and }a_1=s/2)\\
&\text{or }(F=D_n\text{ and }a_1+a_2=s).
\end{aligned}
\end{aligned}
\label{eq:A-J}
\end{equation}

\begin{prop}\label{prop:A-stabilizer}
One has
\[
\begin{aligned}
\widehat\Pi_0(\Lambda)
&=\widehat\Pi(\Lambda)\cap\widehat\Delta_0(\Lambda)
=\{\beta_j:j\in J\},\\
\widehat W_0(\Lambda)&=\langle\sigma_j:j\in J\rangle.
\end{aligned}
\]
Thus the stabilizer is the finite standard parabolic subgroup
generated by the reflections in the walls containing \(x_\Lambda\).
\end{prop}
\begin{pf}
Put \(y=\Lambda+\widehat\rho\). Since \(y\) lies in the closed
fundamental chamber, its stabilizer is generated by the simple
reflections fixing it. Indeed, if \(w\ne1\) fixes \(y\), choose
a left descent \(s_\beta\) of \(w\), with
\(w^{-1}\beta=-\gamma\) and \(\gamma\) positive. Then
\(0\le\langle y,\beta^\vee\rangle
=-\langle y,\gamma^\vee\rangle\le0\), so \(s_\beta\) fixes
\(y\), and \(s_\beta w\) has smaller length. Induction proves the
assertion.

For \(1\le i<n\), the pairing \(\langle y,\beta_i^\vee\rangle\)
is zero precisely when \(a_i=a_{i+1}\); for \(i=n\), it is zero
precisely when \(a_n=0\). The pairing with \(\beta_0^\vee\) is zero
precisely when \(a_1=s/2\) for \(F=\ZZ^n\), or
\(a_1+a_2=s\) for \(F=D_n\). These are exactly the conditions
in \eqref{eq:A-J}, giving
\(\widehat W_0(\Lambda)=\langle\sigma_j:j\in J\rangle\).
Every positive coroot is a nonnegative linear combination of the
simple coroots. Since their pairings with \(y\) are nonnegative,
the coroots with zero pairing are precisely those in the subsystem
generated by \(\{\beta_j^\vee:j\in J\}\). Hence
\(\widehat\Pi_0(\Lambda)=\{\beta_j:j\in J\}\).
Finally, the stabilizer of a point in a discrete Euclidean reflection
group is finite, so \(\widehat W_0(\Lambda)\) is finite.
\end{pf}

Let \(\nu\in\mathcal S^X_{p,n}\).
Since the coefficient of \(E^\rho(\nu)\) in \(M_{p,D_p}\)
is nonzero, there exist \(r,d\in\ZZ_{\ge1}\),
\(\alpha\in\Delta\), and \(w'\in W\) such that
\begin{equation}
rd=D_p,\qquad
2\left(\frac r2\|\alpha\|^2-(\rho|\alpha)\right)=pd,
\qquad
w'(\rho-r\alpha)=\rho+\nu.
\label{eq:A-det}
\end{equation}
Hence
\begin{equation}
\|\rho+\nu\|^2-\|\rho\|^2=pD_p.
\label{eq:A-norm}
\end{equation}
The root \(qd\delta+\alpha\) is integral by
\eqref{eq:A-pairing} and \eqref{eq:A-det}. Reflection in it sends the
finite part of \(\kappa\Lambda_0+\widehat\rho\) from \(\rho\) to
\(\rho-r\alpha\) and changes the coefficient of \(\delta\) by
\(-qD_p\). Together with \(w'\) in \eqref{eq:A-det}, this shows that
\(\kappa\Lambda_0\) and \(\lambda_\nu\) belong to the same integral
Weyl group orbit under the dot action.
In particular, the finite part \(\rho+\nu\) of
\(\lambda_\nu+\widehat\rho\) belongs to
\(W(\rho+sF)\). Hence \eqref{eq:A-common-representative} gives
\(w_{\rho+\nu}^{-1}r(\rho+\nu)=w_\rho^{-1}r(\rho)=x_\Lambda\).

By \eqref{eq:A-orbit-systems}, the three weights have the same
integral root system and integral Weyl group:
\[
\begin{aligned}
\widehat\Delta(\Lambda)&=\widehat\Delta(\kappa\Lambda_0)
=\widehat\Delta(\lambda_\nu),\\
\widehat W(\Lambda)&=\widehat W(\kappa\Lambda_0)
=\widehat W(\lambda_\nu).
\end{aligned}
\]
Their stabilizers are conjugate. Proposition~\ref{prop:A-stabilizer}
therefore shows that all the stabilizers are finite parabolic subgroups. Notice that  the
stabilizers of \(\kappa\Lambda_0\) and \(\lambda_\nu\) need not be
standard for the simple system \(\widehat\Pi(\Lambda)\).

\begin{theorem}\label{thm:A-cosets}
Let \(\nu\in\mathcal S^X_{p,n}\) and
\(\lambda_\nu=\kappa\Lambda_0-qD_p\delta+\nu\). Define
\begin{equation}
\begin{aligned}
\widetilde y&=t_{-c\eta(\rho)}w_\rho,\\
\widetilde z&=t_{-c\eta(\rho+\nu)}w_{\rho+\nu},\\
\widetilde z_0
&=t_{c\left(w_{\rho+\nu}w_\rho^{-1}\eta(\rho)
-\eta(\rho+\nu)\right)}w_{\rho+\nu}w_\rho^{-1}.
\end{aligned}
\label{eq:A-coset-elements}
\end{equation}
These elements belong to \(\widehat W(\kappa\Lambda_0)\), and
\[
\widetilde y\circ\Lambda=\kappa\Lambda_0,\qquad
\widetilde z\circ\Lambda=\lambda_\nu.
\]
The elements \(\widetilde y\) and \(\widetilde z\) are the unique
longest elements of the left cosets
\(\widetilde y\widehat W_0(\Lambda)\) and
\(\widetilde z\widehat W_0(\Lambda)\), respectively. Moreover,
\[
\widetilde z_0=\widetilde z\widetilde y^{-1},\qquad
\widetilde z_0\circ(\kappa\Lambda_0)=\lambda_\nu,
\]
and \(\widetilde z_0\) is the unique shortest element of the left
coset \(\widetilde z_0\widehat W_0(\kappa\Lambda_0)\).
\end{theorem}

\begin{pf}
Since \(\eta(\rho),\eta(\rho+\nu)\in F\) and \(F\) is
\(W\)-stable, all three translation vectors in
\eqref{eq:A-coset-elements} belong to \(cF\). Thus the elements
belong to \(T_{cF}\rtimes W=\widehat W(\kappa\Lambda_0)\).
Equation~\eqref{eq:A-semidirect} gives
\(\widetilde z_0=\widetilde z\widetilde y^{-1}\).

Since \(\widetilde y=h_\rho\), the calculation following
\eqref{eq:A-Lambda} gives
\(\widetilde y\circ\Lambda=\kappa\Lambda_0\).
For \(\widetilde z=h_{\rho+\nu}\), we use
\eqref{eq:A-common-representative} with \(x=\rho+\nu\). The finite part
of \(h_x(\Lambda+\widehat\rho)\) is
\[
w_xx_\Lambda-(\kappa+\mathbf{h}^\vee)c\eta(x)
=r(x)-s\eta(x)=x.
\]
Preservation of the bilinear form determines its \(\delta\)-coefficient
as
\[
\frac{\|\rho\|^2-\|\rho+\nu\|^2}{2(\kappa+\mathbf{h}^\vee)}
=-qD_p
\]
by \eqref{eq:A-norm}. Hence
\(\widetilde z\circ\Lambda=\lambda_\nu\). Since
\(\widetilde z_0=\widetilde z\widetilde y^{-1}\) and
\(\widetilde y\circ\Lambda=\kappa\Lambda_0\), we obtain
\(\widetilde z_0\circ(\kappa\Lambda_0)=\lambda_\nu\).

Lemma~\ref{lem:A-root-signs}, applied to \(x=\rho\) and \(x=\rho+\nu\),
gives
\[
\widetilde y^{-1}\widehat\Delta_0^+(\kappa\Lambda_0)
\subset\widehat\Delta^-,\qquad
\widetilde z^{-1}\widehat\Delta_0^+(\lambda_\nu)
\subset\widehat\Delta^-.
\]
By \eqref{eq:A-inverse-criterion}, \(\widetilde y\) and
\(\widetilde z\) are the longest elements of
\(\widetilde y\widehat W_0(\Lambda)\) and
\(\widetilde z\widehat W_0(\Lambda)\), respectively. Finally, if
\(\alpha\in\widehat\Delta_0^+(\kappa\Lambda_0)\), then
\(\beta=-\widetilde y^{-1}\alpha\in\widehat\Delta_0^+(\Lambda)\)
and
\[
\widetilde z_0\alpha=-\widetilde z\beta\in\widehat\Delta^+.
\]
Lemma~\ref{lem:A-coset-root-criterion} therefore shows that
\(\widetilde z_0\) is the shortest element of
\(\widetilde z_0\widehat W_0(\kappa\Lambda_0)\).
\end{pf}

\vskip 0,3cm
\noindent
{\bf FUNDING}
This work is supported by the National Natural Science Foundation of China
[grant number 12571028].

\vskip 0,3cm
\noindent
{\bf DECLARATION} No AI is  used in this article.

\end{document}